\documentclass[a4paper,10pt]{article}
\usepackage[toc,page]{appendix}
\usepackage{amsfonts}
\usepackage{mathrsfs}
\usepackage{amsmath}
\usepackage{amssymb}
\usepackage{verbatim}
\usepackage{epsfig}
\usepackage{cite}
\usepackage{amstext}
\usepackage{tabularx,booktabs}
\usepackage{graphicx}
\usepackage{subfig}
\usepackage{float}
\usepackage{color}
\usepackage{hyperref}
\usepackage{esint}
\usepackage{authblk}
\usepackage[toc,page]{appendix}
\usepackage{blkarray}
\usepackage{multirow}
\usepackage[intoc]{nomencl}
\makenomenclature
\hypersetup{
    colorlinks,
    citecolor=black,
    filecolor=black,
    linkcolor=black,
    urlcolor=black
}
\numberwithin{equation}{section}
\numberwithin{figure}{section}

\newcolumntype{C}{>{\centering\arraybackslash}X} 
\usepackage{lipsum}
\usepackage{listings}

\definecolor{mGreen}{rgb}{0,0.6,0}
\definecolor{mGray}{rgb}{0.5,0.5,0.5}
\definecolor{mPurple}{rgb}{0.58,0,0.82}
\definecolor{backgroundColour}{rgb}{0.95,0.95,0.92}

\lstdefinestyle{CStyle}{
    backgroundcolor=\color{backgroundColour},   
    commentstyle=\color{mGreen},
    keywordstyle=\color{magenta},
    numberstyle=\tiny\color{mGray},
    stringstyle=\color{mPurple},
    basicstyle=\footnotesize,
    breakatwhitespace=false,         
    breaklines=true,                 
    captionpos=b,                    
    keepspaces=true,                 
    numbers=left,                    
    numbersep=5pt,                  
    showspaces=false,                
    showstringspaces=false,
    showtabs=false,                  
    tabsize=2,
    language=C
}
\usepackage[utf8]{inputenc}
    \usepackage[T1]{fontenc}
    \usepackage{lmodern}
    \usepackage{mathtools}

\begin{document}
\title{On the EOS formulation for  light scattering.\\
  Stability, Singularity and Parallelization}
\author{Aihua Lin}
\author{Per Kristen Jakobsen}
\affil{Department of Mathematics and Statistics, UIT the Arctic University of Norway, 9019 Troms\o, Norway}
\date{\vspace{-5ex}}
\renewcommand\Authands{ and }
\maketitle

\begin{abstract}

\end{abstract}
In this paper we discuss some of the mathematical and numerical issues that have to be addressed when  calculating wave scattering using the EOS approach. The discussion is framed in context of light scattering by objects whose optical response can be of a nonlinear and/or inhomogeneous nature. The discussions address two issues that, more likely than not, will be part of any investigation of wave scattering using the EOS approach.

\section{\bigskip Introduction}
A new hybrid numerical approach for solving  linear and nonlinear scattering problems, the Ewald Oseen Scattering(EOS) formulation, has recently been introduced and applied to the cases of 1D transient wave scattering \cite{Aihua} and 3D light scattering \cite{Aihua2}.
The approach combines a domain-based method and a boundary integral representation in such a  way that the wave fields inside the scattering objects are updated in time using the domain-based method, while the integral representation is used to update the boundary values of the fields, which are required by the  inside domain-based method.  In such a way, for the numerical implementations, no numerical grids outside the scattering objects are needed. 
This greatly reduces the computational complexity and cost compared to fully domain based methods like the Finite Difference Time Domain(FDTD) method or the Finite Element Methods. The method can handle inhomogeneous and/or nonlinear optical response,  and include the time dependent Boundary Element Method(TBEM), as a special case.

For the case of  1D transient wave scattering \cite{Aihua}, the method   solves the model equations accurately and efficiently, but we don't expect the 1D case to be fully representative for the problems and issues that need to be resolved,  while using the EOS formulation to calculate wave scattering.  We do, however, expect the case of 3D light scattering \cite{Aihua2} to be fairly representative with respect to which problems arise, and also the computational and mathematical severity of these problems. We have seen three types of mathematical and computational issues arise for the case of light scattering which we believe are to be found in any nontrivial application of the EOS formulation to wave scattering.

Firstly, we have the issue of numerical stability. Instabilities in numerical implementations of the EOS formulation can arise from discretization of the domain part of the algorithm but also from discretization of the boundary update part of the algorithm.  The  numerical instability arising from the boundary part of the algorithm has been noted earlier in the context of transient  light scattering from objects that has a linear homogeneous optical response. For this situation, realized for example in antenna theory, the boundary part of the EOS algorithm can be disconnected from the domain part of the algorithm, which in this case can be discarded.  The  EOS formulation becomes a pure boundary update algorithm which is solving a set  integro-differential equations located on the boundary of the scattering objects. These integro-differential equations, which are the defining equations for TBEM,  are subject to an instability that, in many common situations,  strikes at late times.
This late time instability is a major nuisance, and  has prevented TBEM  from being more widely applied than it is  today. The sources of these instabilities are not yet  fully understood, but we believe that our investigation of light scattering using the EOS approach, gives some new  insight into the origin of these instabilities. 

  Even without a true understanding of the underlying causes of the late time instability, efforts have been made and several techniques have,  over  the last several decades, been developed with the goal of improving  the stabilities of the numerical schemes designed to solve the  integro-differential equations underlying  TBEM.

  Broadly speaking, there are two different directions that has been  pursued. One direction is to delay or remove the late time instability by applying increasing accurate spatial integration schemes \cite{Weile, Weile2,  Weile3, Weile4, Shanker, Walker, Walker2}. For instance Danile. S. Weiler and his co-authors have published a series of articles focused on illustrating the dependence of the stability on the different numerical integration schemes \cite{Weile, Weile2, Weile3, Weile4}.
The other direction is aimed at designing more stable time discretization schemes.  M. J. Bluck and his co-authors developed a stable, but implicit numerical method, \cite{Walker, Walker2} for the integro-differential equations underlying TBEM,  for the case when the magnetic response is the dominating one. These are the so called magnetic field integral equations. 
Some authors have reported some success in mitigating the instability by both making better approximations to the integrals and also applying improved algorithms for the time derivatives\cite{Zhao, Huang}. 

Our work has not been directly aimed at contributing to this discussion, but, as already noted above,  the integro-differential equations discussed by these authors can be seen as a special case of our general EOS approach, and we therefore believe that the insights we have gained on how this long time instability depend on the different pieces of the EOS algorithm, in particular how it depends on the material parameters describing the optical response of the scattering object, do have some relevance to the discussion described above.

Secondly, there is the issue of the singular integrals that appear when the integral part of the EOS algorithm is discretized. This issue is very much present in BEM and in TBEM \cite{singular1,singular2,singular3,singular4}, but they are more prevalent and severe for the EOS formulation, where we have to tackle both surface integrals and volume integrals.  We believe that the type of singular integrals, and how to treat them for the case of light scattering, are fairly representative for the level of complexity one will encounter, while applying the EOS approach to wave scattering problems. For this reason we find it appropriate to include a section in this paper, where we discuss relevant types of integrals, and how to treat them.

Thirdly, the fundamental equations underlying both the TBEM and our more general EOS approach to transient wave scattering, are retarded in time. This retardation is unavoidable since their underlying equations can only be derived using space-time Green's functions. Thus the solutions at a certain time depend on a values of the solutions from a potentially very long previous interval of time. Computationally this means that the method can be very demanding with respect to memory, and it also means that the updating of the boundary values of the fields, which is done by the boundary part of the EOS algorithm, can be very costly. Parallel  processing, either using a computational cluster or a shared memory machine can take on these computational tasks. However, whenever large scale parallel processing is needed, the issue of appropriate partitioning of the problem and load balancing inevitably comes into play. In our work the EOS algorithm was implemented on a large cluster, but we will not in this paper report on any of the parallel issues that our EOS approach for light scattering gave rise to. These kind of considerations, which are important in practical terms, but typically have fairly low generality, are somewhat distinct from the mathematical and numerical issues that are the focus of the current paper, and will therefore be reported elsewhere at a later time.

However,  the high memory requirement of the EOS approach to light scattering, is something that should be addressed at this point.  On the one hand, the EOS approach represents a large, potentially very large,  reduction in memory use, as compared to fully domain based methods, since only the surface and inside of the scattering objects has to be discretized.  On the other hand, because of the retardation, there is a large, potentially very large increase in memory use compared to the memory usage needed by the domain part of the algorithm. It is appropriate to ask if anything has been gained with respect to memory usage compared to a fully domain based method like the FDTD method? We don't, as of yet, know the answer to this question, and the answer is almost certainly  not going to be a simple one. It  will probably depend on the detailed structure of the problems like the nature of the source, the number, shape and distribution of scattering objects etc. However, even if the memory usage for purely domain based methods and our EOS approach are roughly the same for many problems of interest, our approach avoid many of the sources of problems that need to be taken into account while using purely domain based methods. These are problems  like stair-casing  at sharp interfaces defining the scattering objects, issues of accuracy, stability and complexity associated with the use of multiple grids in order to accommodate the possibly different geometric shapes of the scattering objects and the need to minimize the reflection from the boundary of the finite computational box. The EOS approach is not subject to any of these problems.

In this paper our effort are aimed towards testing  the EOS formulations of light scattering with respect to implementation complexity and
numerical stability. Thus we illustrate the method by the simplest situation where we have single scattering object in the form of a rectangular box.

 In section \ref{stabilities3d}  we analyze the numerical stability of our EOS scheme for light scattering by using eigenvalues of the  matrix defining the linearized version of the scheme exactly  like for the case of 1D wave scattering\cite{Aihua}. We find, just like for the 1D case, that the internal numerical scheme, Lax-Wendroff for our case determines a stability interval for the time step. In the 1D case, the stability interval of the EOS formulation is purely determined by the internal numerical scheme. However for the 3D case, there is another lower limit of the stability interval determined by the integral part of the scheme which leads to the situation  where  the lower limit of the stability interval is determined by the integral equations, and the upper limit is determined by the internal numerical scheme.  We find that  the late time instability is highly depended on the features of the scattering materials and specifically, it is directly related to the values of the relative magnetic permeability $\mu_1$ and the relative electric permittivity $\varepsilon_1$.  Using this we prove that, for the relative permeability and permittivity  in a certain range,  the numerical scheme for our EOS formulation of light scattering, works well and is without any late time instabilities.  The late time instability is only observed for high relative electric permittivity or high relative magnetic permeability.
We also observe that the lower limit of the stability interval for the time step is more  sensitive to relative differences in magnetic permeability $\mu_1$  than electric permittivity $\varepsilon_1$ between the inside and outside of the scattering objects.
 
 In section \ref{Asingularity} we present the singular integrals that appear in our EOS formulation for light scattering and the techniques we use to reduce their calculation to a singular core, which we calculate exactly, and a regular part which we calculate numerically.

\section{Stability}\label{stabilities3d}

In this section we discuss instabilities showing up at late times when we discretize the EOS formulation for light scattering. Whether  or not the late time instability show up, depends on the values of the material parameters defining the problem. The overall method is far to complex for an analytical investigation of the stability to be feasible, but using numerical calculation of the eigenvalues of a linearization of the system of difference equations defining the numerical implementation of the EOS formulation, supplemented by running of the full algorithm, we find that the domain part and the boundary part of the algorithm contribute to the instability separately and in different ways. The focus of this section is to disentangle these two contributions to the instability.
For the domain part of the algorithm we use  Lax-Wendroff, which is an explicit method. The discrete grid inside the scattering object must, for the EOS formulation of light scattering, support both discrete versions of the partial derivatives, and also discretizations of the integrals defining the boundary update  part of the algorithm. For this reason the grid is nonuniform close to the boundary. The discretization of the domain part of the algorithm takes the form of a vector iteration
\begin{equation}\label{3dmatrixequation}
Q^{n+1}=MQ^n,
\end{equation}
where $Q$ is a vector containing the components of the electric field and the magnetic field at all points of the grid with a size $6\times N_x \times N_y \times N_z$, where $N_x, $ $N_y$ and $N_z$ are the number of grid points in the  $x,$ $ y$ 
and $z$ directions.
The entries of the matrix $M$ are presented in Appendix \ref{Amatrix}.
In order to get a stable numerical solution, as discussed in \cite{Aihua}, the largest eigenvalues of the matrix $M$ must have a norm smaller than 1. 
For the non-uniform grids and the discretizations in \cite{Aihua2}, we find that the vector iteration (\ref{3dmatrixequation}) is stable if
$$
0.005< \tau <0.48,
$$ 
where $\tau=c_1 \Delta t/\Delta x.$
\begin{figure}[h!]
\centering
\includegraphics[width=8cm]{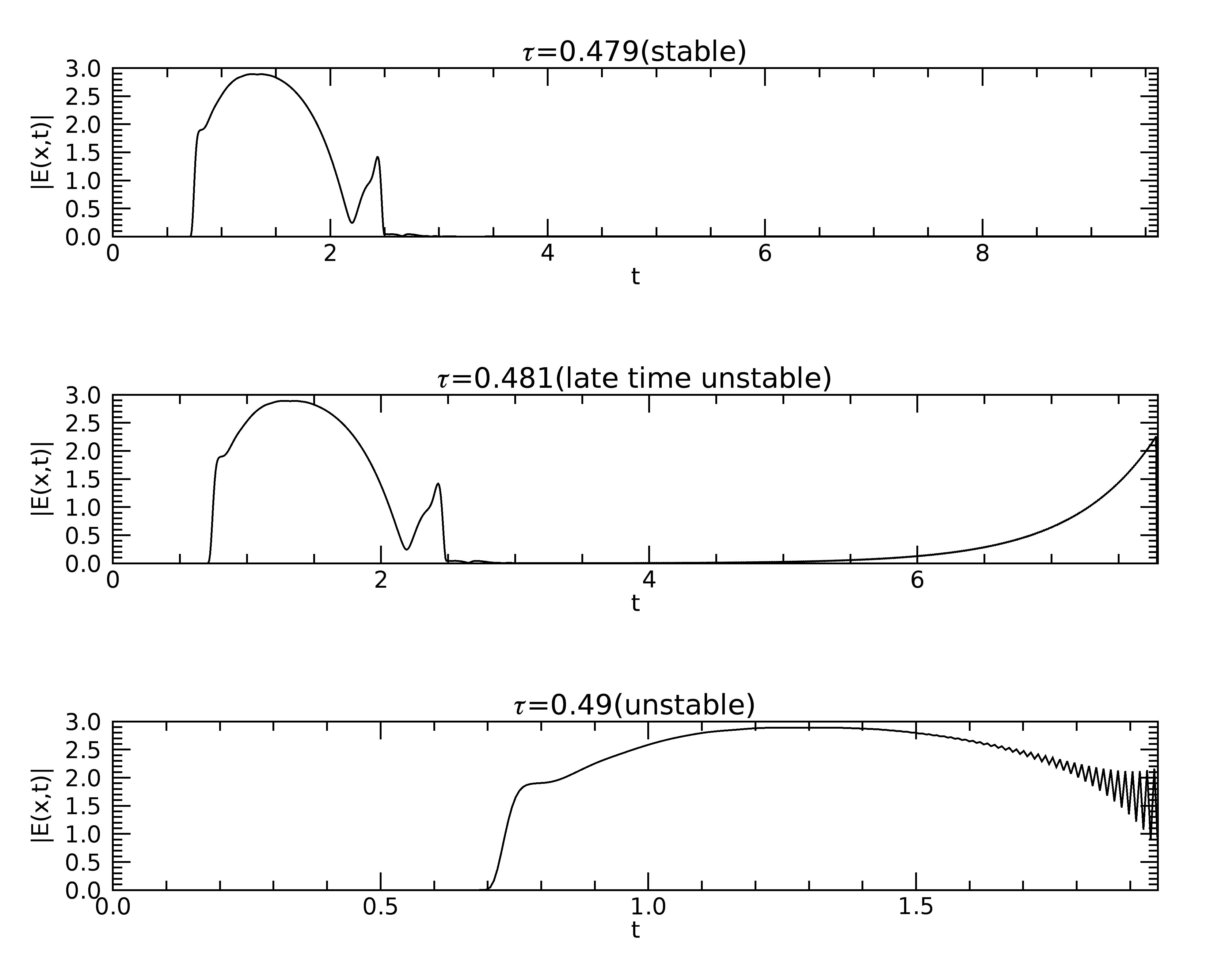}
\caption{Numerical solutions from different values of $ \tau .$ $\mu_1=1.0, \varepsilon_1=1.5, \mu_0=1.0, \varepsilon_0=1.0.$ }
\label{unstablei}
\end{figure}

\noindent Figure \ref{unstablei} illustrates the intensity of the electric field at a specific point inside the object, as a function of time,  for different values of $\tau$.  The instability, which in the TBEM literature is called the late time instability, is illustrated in the second panel of figure \ref{unstablei}. As we mentioned in the introduction in the paper, the term late time instability has been much used  in the community that is focused on time dependent boundary element method. We believe that in their domain of application, like antenna theory, the physical parameters are such that the largest eigenvalue for the iteration  is always only slightly bigger than 1, like it is in panel two of figure \ref{unstablei} . That's why the instability always shows up at late times. In panel three of the figure we are deeper into the unstable domain for $\tau$, and the larges eigenvalue is now so large that  it destroys the whole calculation. The late time instability has thus been transformed into an early time instability.  Note that the outside source in figure \ref{unstablei} is the same as in \cite{Aihua2}.

\begin{figure}[h!]
\centering
\includegraphics[width=8cm]{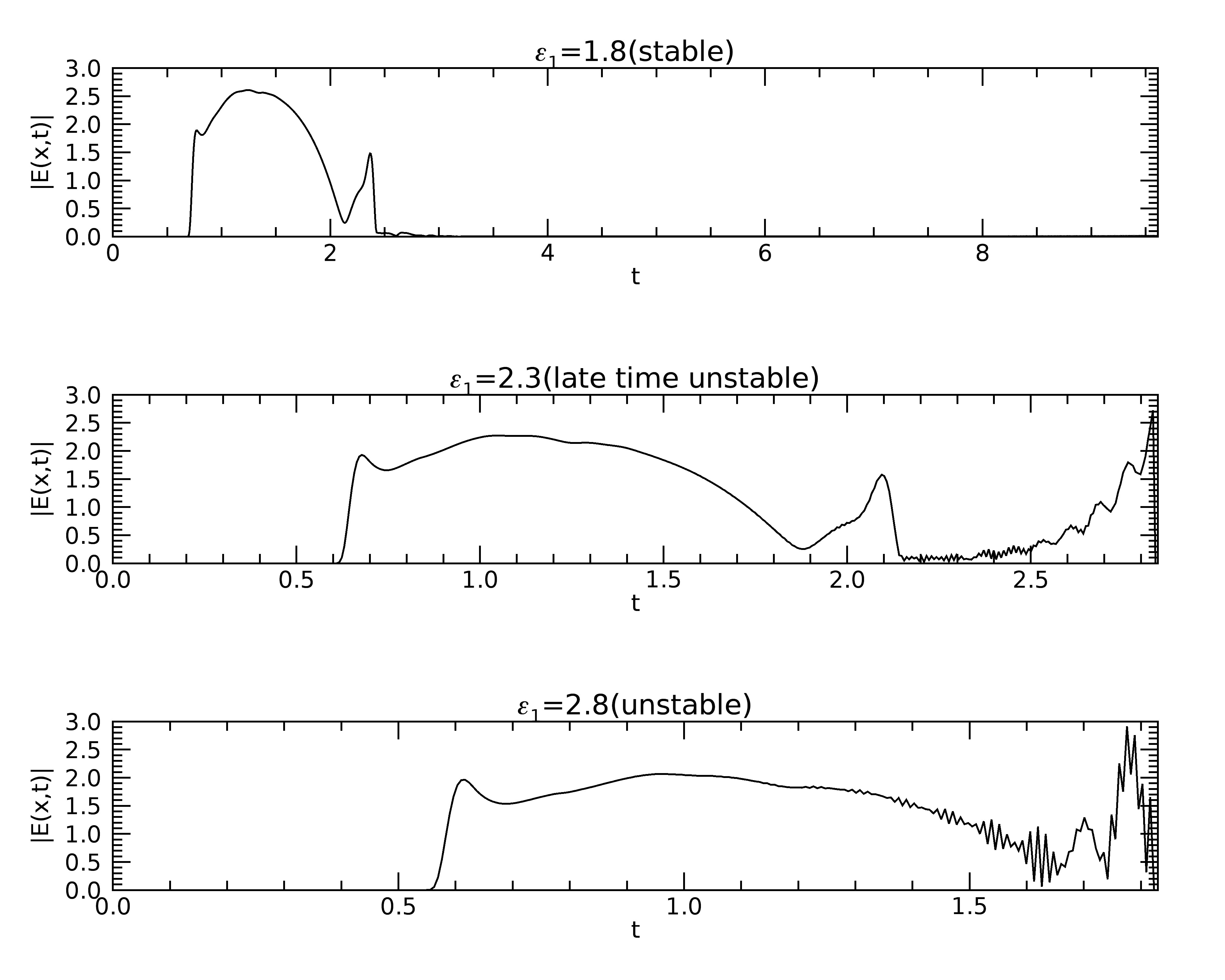}
\caption{Numerical solutions from different values of $ \varepsilon_1 .$ $\tau=0.45,$ $\mu_1=1.0, \mu_0=1.0, \varepsilon_0=1.0.$ }
\label{unstables}
\end{figure}

\begin{figure}[h!]
\centering
\includegraphics[width=8cm]{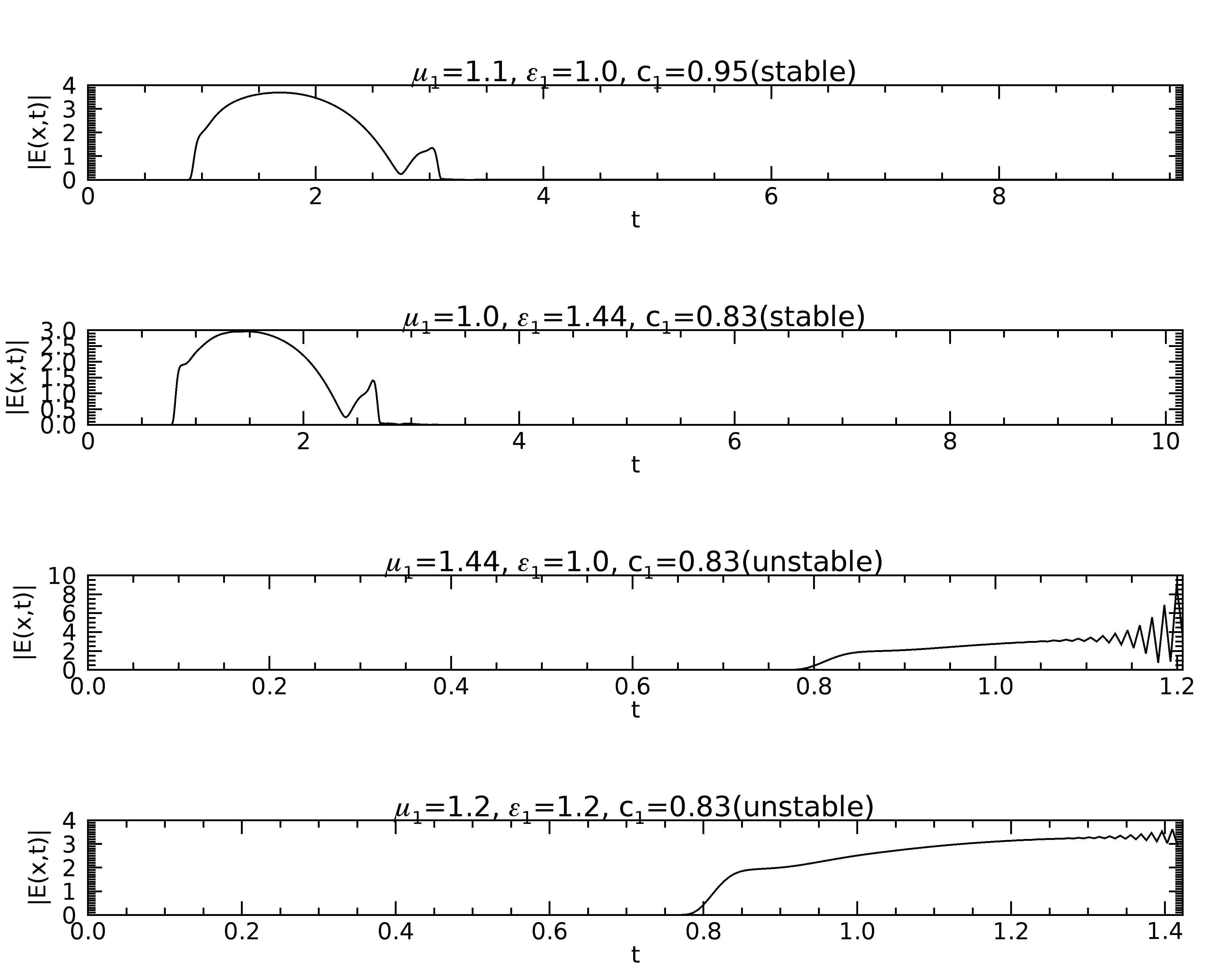}
\caption{Numerical solutions from different values of $ \mu_1 .$ $\tau=0.45,$ $\mu_0=1.0, \varepsilon_0=1.0.$}
\label{ucompare}
\end{figure}

In our numerical experiments, we found that the stable range of the EOS formulations is not only restricted by the eigenvalues of the matrix $M$, but is also restricted by the boundary integral identities through the relative electric permittivity $\varepsilon_1$ and the relative magnetic permeability $\mu_1$. Figure \ref{unstables} shows how the stability depends on the values of $\varepsilon_1$, and figure \ref{ucompare} shows how it depends on the values of $\mu_1$. Together, they tell us that increasing the electric permittivity or the magnetic permeability narrows the stable range. 
  Figure \ref{ucompare} also tells us that $\mu_1$ and $\varepsilon_1$ don't affect the stability of the full scheme in the same way. It seems that the method is more sensitive to $\mu_1$ than $\varepsilon_1.$
After a series of numerical experiments, our conclusion is that, for an explicit numerical method like the one we are using, the lower limit of the stable range of the EOS formulation is restricted by the electric permittivity $\varepsilon_1$ and the magnetic permeability $\mu_1$ while the upper limit of the stable range is determined by the inside domain-based method. This conjecture is verified by the following two tests.

\subsection{Instabilities coming from the domain-based method}
 For the first test we consider a homogeneous model without current and charge inside the object  which implies $\mu_1=\mu_0,$ $\varepsilon_1=\varepsilon_0,$ ${\bf J}_1=0$ and $\rho_1=0.$ Under these assumptions, the electric field and the magnetic field are continuous across the surfaces,
\begin{equation*}
\begin{split}
{\bf E}_-&={\bf E}_+,\\
{\bf B}_-&={\bf B}_+,
\end{split}
\end{equation*}
where ${\bf E}_\pm$and ${\bf B}_\pm$ are the integral representations of the solutions on the surface by taking the limit from the inside and the outside of the object respectively.
The electric field inside the object can be calculated by the outside sources directly
\begin{equation}\label{htestequation2}
\begin{split}
{\bf E}_1({\bf x},t)=-\partial_t\frac{\mu_0}{4\pi}\int_{V_0}\,\mathrm{d}V'\frac{{\bf J}_0({\bf x'},T)}{|\bf x'-\bf x|}-\nabla\frac{1}{4\pi\varepsilon_0}\int_{V_0}\,\mathrm{d}V'\frac{\rho_0({\bf x'},T)}{|{\bf x'}-{\bf x}|},
\end{split}
\end{equation}
where ${\bf x}\in V_1.$ (\ref{htestequation2}) expresses the exact solution for the inside fields.
Also from  \cite{Aihua2} we have the boundary integral identity
\begin{equation}
\begin{split}\label{htestb}
{\bf E}_+({\bf x},t)=-\partial_t\frac{\mu_0}{4\pi}\int_{V_0}\,\mathrm{d}V'\frac{{\bf J}_0({\bf x'},T)}{|\bf x'-\bf x|}-\nabla\frac{1}{4\pi\varepsilon_0}\int_{V_0}\,\mathrm{d}V'\frac{\rho_0({\bf x'},T)}{|{\bf x'}-{\bf x}|},
\end{split}
\end{equation}
and 
\begin{equation}\label{htestbb}
{\bf B}_+({\bf x},t)=\nabla\times\frac{\mu_0}{4\pi}\int_{V_0}\,\mathrm{d}V'\frac{{\bf J}_0({\bf x'},T)}{|{\bf x'}-{\bf x}|},
\end{equation}
for ${\bf x} \in S$, ${\bf E}_+({\bf x},t)$ and ${\bf B}_+({\bf x},t)$ represent the limits by letting ${\bf x}$ approach the surface from the inside of the scattering object. 
On the other hand,  \cite{Aihua2} gives the integral representations for the inside domain by 
\begin{equation}\label{htestequation1}
\begin{split}
{\bf E}_1({\bf x},t)&=\partial_t[\frac{1}{4\pi}\int_{S}\,\mathrm{d}S^{'}\{\frac{1}{c_1|{\bf x'}-{\bf x}|} ({{\bf n}'}\times {\bf E}_+({\bf x'},T))\times \nabla'|{\bf x'}-{\bf x}|\\
&+\frac{1}{c_1|{\bf x'}-{\bf x}|} ({{\bf n}'}\cdot {\bf E}_+({\bf x'},T))\nabla'|{\bf x'}-{\bf x}|+\frac{1}{|{\bf x'}-{\bf x}|} {{\bf n}'}\times {\bf B}_+({\bf x'},T)  \}]\\
&-\frac{1}{4\pi}\int_{S}\,\mathrm{d}S^{'}\{ ({{\bf n}'}\times {\bf E}_+({\bf x'},T))\times \nabla'\frac{1}{|{\bf x'}-{\bf x}|}\\
&+({{\bf n}'}\cdot {\bf E}_+({\bf x'},T))\nabla'\frac{1}{|{\bf x'}-{\bf x}|}\}.
\end{split}
\end{equation}
\begin{figure}[h]
\centering
\includegraphics[width=8cm]{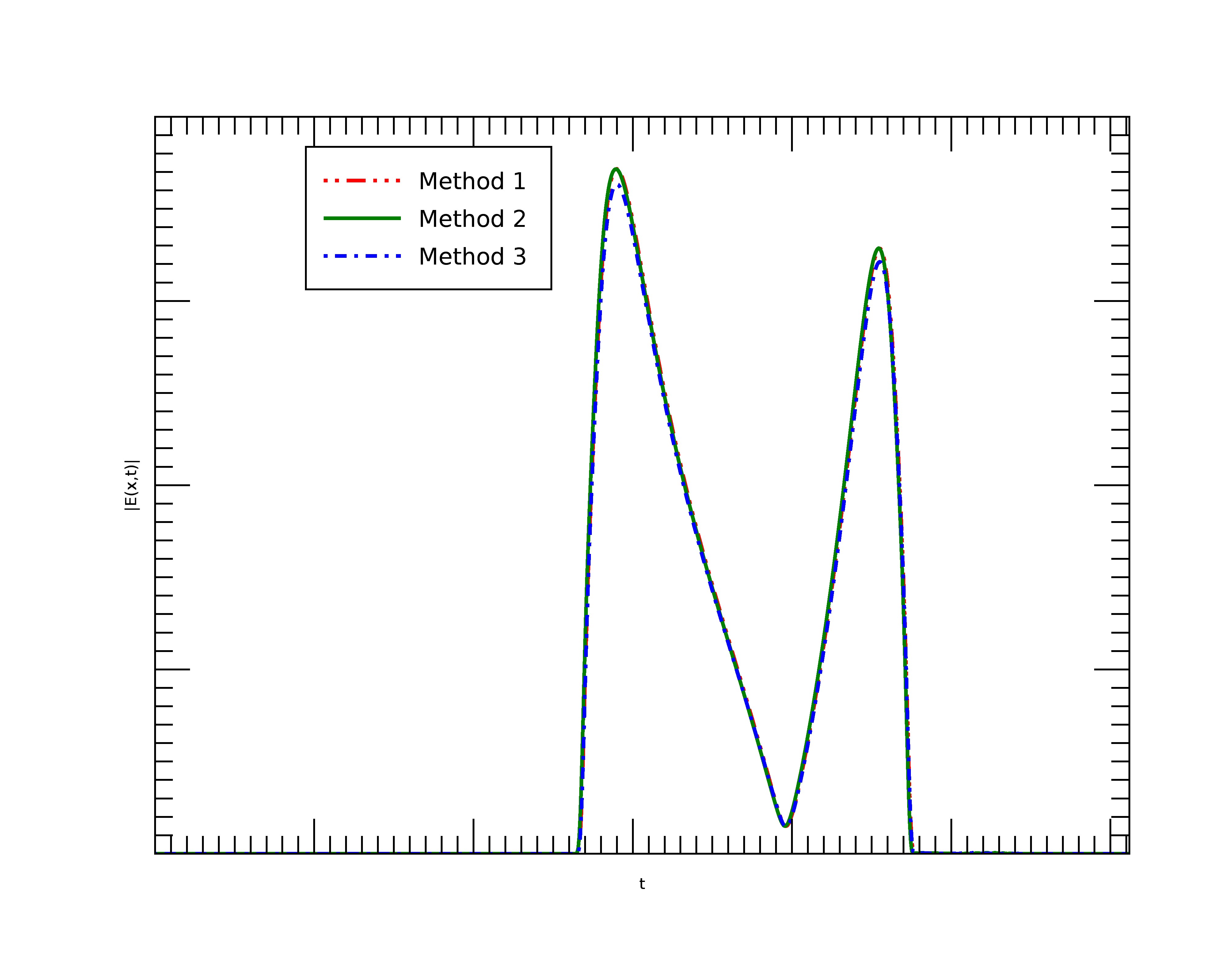}
\caption{Comparison of  the intensity of the electric field inside the object at a specific point calculated by three methods.   $t_0=1.5,$  $x_0=-2.0,$  $y_0=0.0,$  $z_0=0.0,$  $\tau=0.45,$  $\mu_1=1.0,$ $ \varepsilon_1=1.0,$  $\mu_0=1.0, $ $\varepsilon_0=1.0.$ }
\label{htest1}
\end{figure}
Thus the solution for the domain inside the scattering object can now be calculated in three ways. The first is the exact solution expressed by (\ref{htestequation2}), the second, Method 2,  is the Lax-Wendroff method supplied by the exact boundary values (\ref{htestb}) and (\ref{htestbb}) , and the third, Method 3,  is to calculate the solution using formula (\ref{htestequation1}) which expresses the field values inside the scattering object in terms of the values of the fields on the boundary.  Note that Method 3 uses the same surface integral expressions as the one that form the boundary  part of the full implementation of our EOS formulation of light scattering. Thus, instabilities in the full algorithm originating from the boundary part of the algorithm, should appear as instability in Method 3.  

Figure \ref{htest1} compare the solutions calculated  in these three ways,  where $\mu_1, \nu_1$ and $\tau$ have been fixed in the stable range.  Both Method 2 and Method 3  are stable and give solutions that agree with the exact solution to high accuracy.
\begin{figure}[h!]
\centering
\includegraphics[width=8cm]{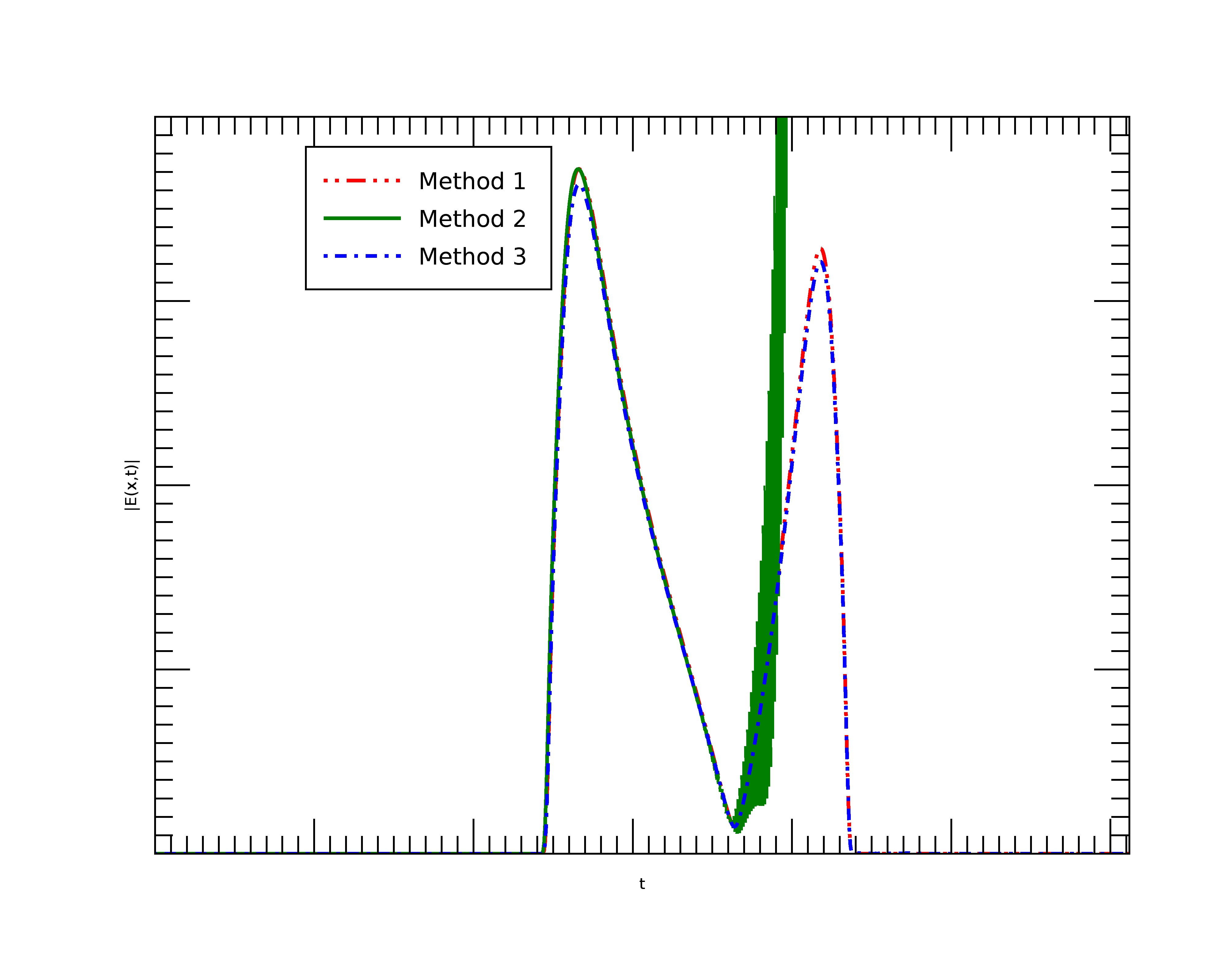}
\caption{Comparison of  the intensity of the electric field inside the object at a specific point calculated by three methods.   $t_0=1.5,$  $x_0=-2.0,$  $y_0=0.0,$  $z_0=0.0,$  $\tau=0.49,$  $\mu_1=1.0,$ $ \varepsilon_1=1.0,$  $\mu_0=1.0, $ $\varepsilon_0=1.0.$}
\label{htest2}
\end{figure}
 In Figure \ref{htest2}$, \tau$ has been set to be 0.49, and is thus is larger than  the upper limit of the stable range. The figure shows that Method 2 is now unstable but Method 3 is still stable and equal to the exact solution to high accuracy. The outside source in figure \ref{htest1} and figure \ref{htest2}  is as same as in \cite{Aihua2} and the values of the parameters are shown under the figure. 
\subsection{Instabilities coming from the boundary integral identities}
In order to investigate the dependence of the stability on $\mu_1$ and $\varepsilon_1,$  we set up a test based on the use of artificial sources as  in \cite{Aihua2}. The idea is to chose functional forms for an electromagnetic field, and then calculate the sources, charge density and current density,  needed for making the chosen fields solutions to Maxwell's equations driven by the calculated sources

We now calculate the electromagnetic field inside the scattering object in two different ways.  In Method 1 we use the discretization of the EOS formulation developed in \cite{Aihua2}, which combines the Lax-Wendroff method for the domain part of the algorithm and our discretization of the integral representations of the boundary fields for the boundary part of the algorithm.   Method 2  is to calculate the inside field values by only using the Lax-Wendroff method supplemented by the exact boundary values of the electromagnetic field which are the ones we chose while setting up the artificial sources.  
\begin{figure}[h!]
\centering
\includegraphics[width=8cm]{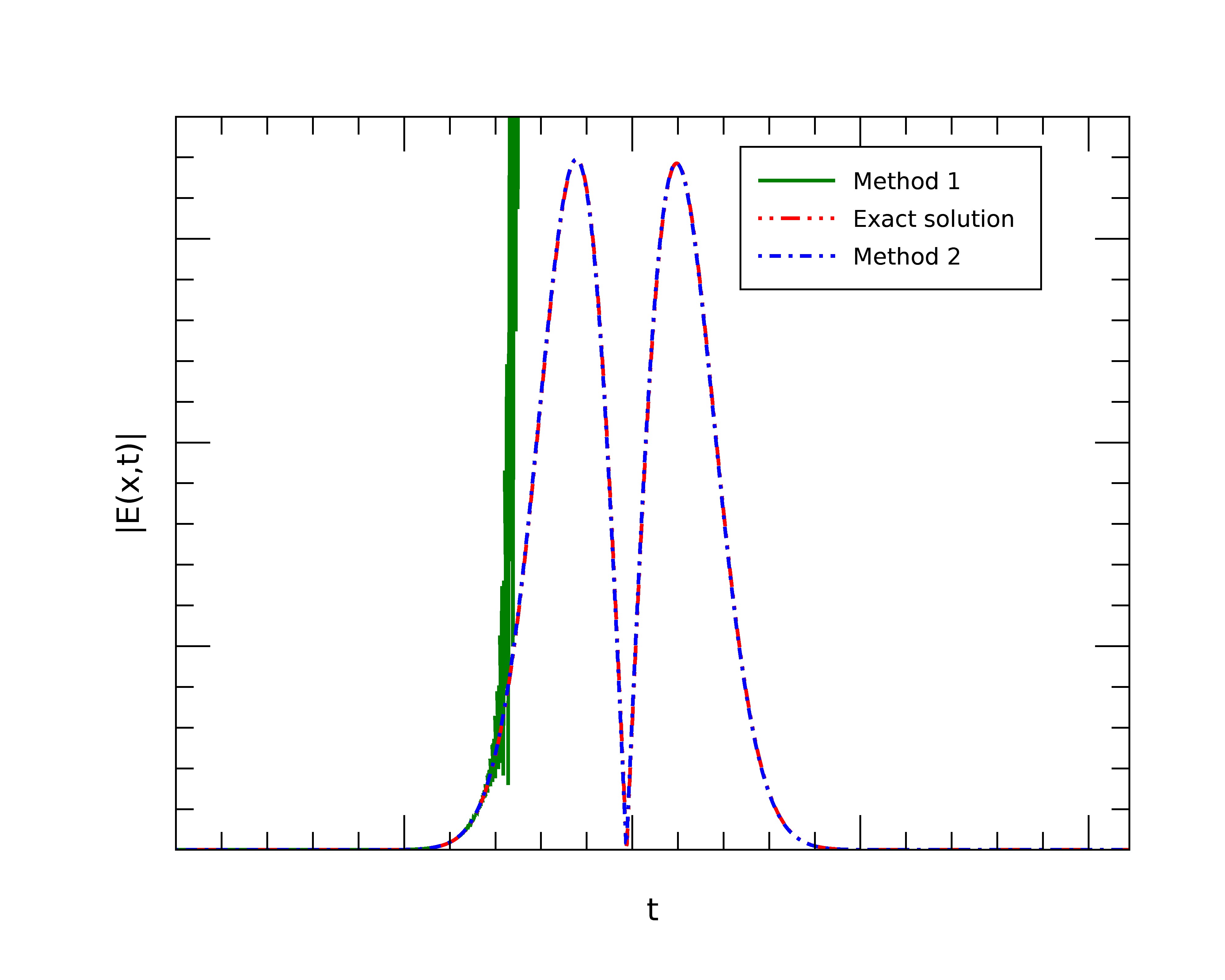}
\caption{Comparison of  the intensity of the electric field inside the object at a specific point between the exact solution and the numerical results calculated by two methods.  $\tau=0.45, $ $\mu_1=1.0,$ $ \varepsilon_1=2.5,$  $\mu_0=1.0, $ $\varepsilon_0=1.0.$}
\label{uhtest}
\end{figure}
Figure \ref{uhtest} is the numerical result where the upper limit of the stable range is kept  while the values of $\mu_1$ and $\varepsilon_1$ have been chosen to break the lower limit of the stable range of the EOS formulations. It shows that even though the lower limit of the stable range has been broken, Method 2, which only involves the Lax-Wendroff method works perfectly. \ref{htest2} and \ref{uhtest} tell us that the changing of the lower limit does not effect the stability of the Lax-Wendroff method and the changing of upper limit does not effect the stability of the surface integrals. For a general application where the source is located outside the object and there are current density and electric density inside the scattering object, the EOS formulations does have a range for a stable numerical implementation. The upper limit of the range is determined by the Lax-Wendroff method due to the non-uniform grids and the lower limit is determined by the changing $\mu_1$ and $\varepsilon_1.$ The setting up of the artificial sources and the values of the parameters in figure \ref{uhtest} are the same as the artificial sources in \cite{Aihua2}. From figure \ref{htest2} and figure \ref{uhtest}, we can also see that before the instabilities show up, both the EOS formulations and the Lax-Wendroff method solve the equations accurately.

\section{Calculations of the singular integrals}\label{Asingularity}
In this section we introduce a technique to accurately calculate  integrals with singularities which can be applied for both the singular volume integrals and the singular surface integrals occurring in the EOS formulations of the 3D Maxwell's equations.  Here we illustrate the technique by calculating one type of singular volume  integral
\begin{equation}\label{f1}
f_1=\iiint_{V_{i,j,k}} \frac{1}{|{\bf x}'-{\bf x}_p|}\,\mathrm{d}V=\iiint_{V_{i,j,k}} \frac{1}{r}\,\mathrm{d}V,
\end{equation}
where the integral domain  $V_{i,j,k}$ is adjacent to the surfaces of the scattering object and given by
\begin{equation*}
V_{i,j,k}=[x_a, x_a+\Delta x]\times [y_j-\frac{\Delta y}{2}, y_j+\frac{\Delta y}{2}] \times [z_k-\frac{\Delta z}{2}, z_k+\frac{\Delta z}{2}],
\end{equation*}
with surfaces $S_m, $ $ m=1,2,\cdots,6$.   Here, $\Delta x,$ $\Delta y$ and $\Delta z$ are the grid parameters  in $x,$ $y$ and $z$ directions respectively.

\begin{figure}[h!]
\centering
\includegraphics[width=6cm]{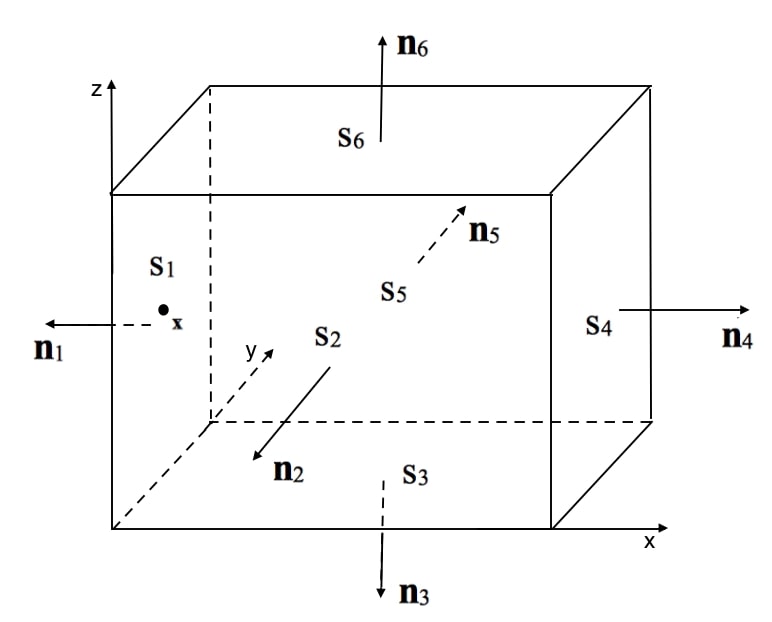}
\caption{The integral  domain of the singular integral }
\label{box}
\end{figure}

The point  ${\bf x}_p$
$${\bf x}_p=(x_a, y_j, z_k),$$
 is centered on one of the  surfaces of the scattering object. 
 The geometry is illustrated in figure \ref{box}, where ${\bf n}_m$ is the unit normal vector on surface $S_m$ pointing out of $V_{i,j,k}.$ 

The components of the integration variable in (\ref{f1}) are given by
$$
{\bf x'}=(x', y', z'),
$$
and let us introduce the quantity
$$
{\bf r}={\bf x}'-{\bf x}_p,
$$
with $r=|{\bf r}|.$

We want to apply the divergence theorem on (\ref{f1}), and therefore need to find a function $\varphi(r)$ that satisfies 
\begin{equation*}
\nabla \cdot ({\bf r} \varphi(r))=\frac{1}{r},
\end{equation*}
or equivalently
$$3 \varphi(r)+r\varphi'(r)=\frac{1}{r}.$$
Solving the above equation, we get 
\begin{equation*}
\varphi(r)=\frac{1}{2 r}.
\end{equation*}
Because of the singularity on $S_1,$ we can not apply the divergence theorem directly, however we can write $f_1$ as
\begin{equation*}
\begin{split}
f_1= \frac{1}{2 }(\sum_{m=2}^{6} \iint_{S_m} \frac{1}{ r} {\bf r}\cdot {\bf n}_m\,\mathrm{d}S+\lim_{\epsilon \rightarrow 0}\iint_{{S_\epsilon}} \frac{1}{r} {\bf r}\cdot {\bf n}_{\epsilon}\,\mathrm{d}S+\iint_{S_\Omega} \frac{1}{r} {\bf r}\cdot {\bf n}_1 \,\mathrm{d}S),
\end{split}
\end{equation*}
where $S_\epsilon$ is a hemispherical surface of radius $\epsilon$ centered at ${\bf x}_p$ and  $S_\Omega$ is the rest of the surface $S_1$ with a disk of radius $\epsilon$ around ${\bf x}_p$ has been removed. 
${\bf n}_\epsilon$ is the unit normal vector on $S_\epsilon,$ pointing out of $V_{i,j,k}.$ 
${\bf n}_m$ is the unit normal vector on $S_m,$ pointing out of $V_{i,j,k}.$ 

For the integral over $S_\Omega$,  we have 
 $${\bf r}=(0, y'-y_j, z'-z_k)$$ and
$$
{\bf n}_1=(-1,0,0),
$$
thus we get
\begin{equation*}
\iint_{S_\Omega}\frac{1}{ r} {\bf r}\cdot {\bf n}_1\,\mathrm{d}S=0. 
\end{equation*} 
For the integral over $S_\epsilon$, we use the spherical coordinate system, $${\bf r}=\epsilon (\cos \theta \sin \varphi, \sin \theta \sin \varphi, \cos \varphi),$$
and 
$$ {\bf n}_{\epsilon}= (\cos \theta \sin \varphi, \sin \theta \sin \varphi, \cos \varphi),$$
where $\epsilon, \varphi, \theta$ are  respectively the radial distance, polar angle and azimuthal angle,
so that
\begin{equation*}
\begin{split}
\lim_{\epsilon \rightarrow 0} \iint_{S_\epsilon}\frac{1}{ r} {\bf r}\cdot {\bf n}_{\epsilon}\,\mathrm{d}S=&\lim_{\epsilon \rightarrow 0}\frac{1}{\epsilon}\int_{0}^{2\pi}\int_{-\frac{\pi}{2}}^{\frac{\pi}{2}} \epsilon (\cos \theta \sin \varphi, \sin \theta \sin \varphi, \cos \varphi)\\
&\cdot  (\cos \theta \sin \varphi, \sin \theta \sin \varphi, \cos \varphi) \epsilon^2 \sin \varphi \, \mathrm{d} \theta\, \mathrm{d} \varphi\\
&=0.
\end{split}
\end{equation*}
Defining
$$
s_m= \iint_{S_m} \frac{1}{ r} {\bf r}\cdot {\bf n}_m\,\mathrm{d}S,
$$
$f_1$ can be written as
\begin{equation}\label{f1s}
f_1= \frac{1}{2 }\sum_{m=2}^{6} s_m.
\end{equation}
(\ref{f1s}) is not singular any more and can be calculated by 2D Gaussian quadrature. However we will compute $f_1$  by reducing the surface integral into a line integral, which is also the approach we use to calculate the singular surface integrals appearing in the implementation discussed in this paper. 

\begin{figure}[H]
\centering
\includegraphics[width=6cm]{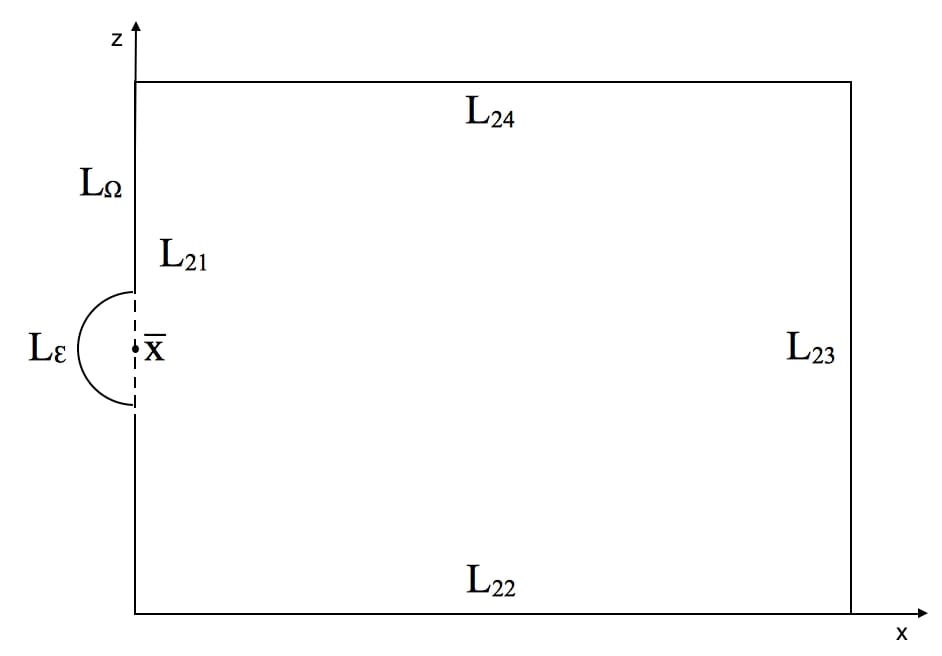}
\caption{Surface $S_2$ }
\label{s2}
\end{figure}

We first consider the integral over $S_2$. The geometry is shown in figure \ref{s2}. 
As show in  figure \ref{s2}, the surface $S_2$  is bounded by the union of four straight lines $L_{2n}$, $n=1,2,3,4.$ 
On this surface we have $${\bf r}=(x'-x_a, \frac{1}{2}\Delta y, z'-z_k)$$ and the unit normal is  $${\bf n}_2=(0,-1,0),$$
so that
\begin{equation*} 
s_2=\frac{1}{2}\Delta y\iint_{S_2} \frac{1}{\sqrt{(x'-x_a)^2+\frac{1}{4}\Delta y^2+(z'-z_k)^2}}\, \mathrm{d}S.
\end{equation*}
The goal is to use the divergence theorem on this surface integral and thereby reduce it to line integrals over the four lines that forms the boundary of $S_2$. We therefore seek a function $\varphi(\bar r)$ that satisfies
 $$\nabla \cdot (\bar {\bf r} \varphi(\bar r))=\frac{1}{\sqrt{\bar r^2+\frac{1}{4}\Delta y^2}},$$
where $\bar {\bf r}=(x'-x_a,z'-z_k)$ and $\bar r=|\bar {\bf r}|.$  This equation can be rewritten in the form
$$2 \varphi(\bar r)+r\varphi'(\bar r)=\frac{1}{\sqrt{\bar r^2+\frac{1}{4}\Delta y^2}}.$$
Solving the above equation we get
\begin{equation*}\label{f1s2}
 \varphi(\bar r)=\frac{\sqrt{\bar r^2+\frac{1}{4}\Delta y^2}}{\bar r^2} .
\end{equation*}

Using the divergence theorem and taking into account of the singularity at $$\bar{ \bf  x}=(x_a, z_k)$$ on $L_{21},$ we get
\begin{equation*} 
\begin{split}
s_2=\frac{1}{2}\Delta y (\sum_{n=2}^4 \int_{L_{2n}} \varphi(\bar r)\bar {\bf r}\cdot \bar {\bf n}_{n}\,\mathrm{d}L
+  \lim_{\epsilon \rightarrow 0}\int_{L_\epsilon} \varphi(\bar r) \bar {\bf r}\cdot \bar {\bf n}_{ \epsilon}\, \mathrm{d}L+  \int_{L_\Omega} \varphi(\bar r) \bar {\bf r}\cdot \bar {\bf n}_{1}\, \mathrm{d}L)
\end{split}
\end{equation*}
where  $L_\epsilon$ is a semicircle with radius $\epsilon$ centered at point ${\bar {\bf  x}}$ and $L_\Omega$ is the rest of $L_{21}$. Here $\bar {\bf n}_{n}$ is the unit normal of $L_{2n}$, pointing out of $S_2,$ and $\bar {\bf n}_{\epsilon}$ is the unit normal of $L_{\epsilon}$, pointing out of $S_2.$ 

For the integral over $L_\Omega,$ 
we have $$\bar {\bf r}=(0, z'-z_k),$$ and $$\bar {\bf n}_1=(-1,0),$$ so that
\begin{equation}\label{h1s21}
\int_{L_\Omega} \frac{\sqrt{\bar r^2+\frac{1}{4}\Delta y^2}}{\bar r^2} (0, z'-z_k)\cdot (-1,0)\, \mathrm{d}L=0.
\end{equation}
For the integral over  $L_\epsilon,$ using the polar coordinates, we have $$\bar {\bf r}=\epsilon (\cos \theta , \sin \theta),$$ and
$$\bar  {\bf n}_\epsilon=-(\cos \theta , \sin \theta),$$
so that 
\begin{equation}\label{h1s22}
\begin{split}
&\lim_{\epsilon \rightarrow 0}\int_{L_\epsilon} \frac{\sqrt{\bar r^2+\frac{1}{4}\Delta y^2}}{\bar r^2} \bar {\bf r}\cdot\bar  {\bf n}_\epsilon\, \mathrm{d}L\\
&=-\lim_{\epsilon \rightarrow 0}\int_{-\frac{\pi}{2}}^{\frac{\pi}{2}} \epsilon (\cos \theta , \sin \theta)\cdot  (\cos \theta , \sin \theta ) \frac{\sqrt{\epsilon^2+\frac{1}{4}\Delta y^2}}{\epsilon^2} \epsilon\, \mathrm{d} \theta\\
&=-\frac{1}{2} \Delta y \pi.
\end{split}
\end{equation}
Summing up (\ref{h1s21}) and (\ref{h1s22}) gives $$l_{21}=-\frac{1}{2} \Delta y \pi.$$
Thus $s_2$ is expressed by
\begin{equation*} 
\begin{split}
s_2=\frac{1}{2}\Delta y\sum_{n=1}^4 l_{2n},
\end{split}
\end{equation*}
where
\begin{equation*}
\begin{split}
 l_{22}&=\frac{1}{2}\Delta z \int_{x_a}^{x_a+\Delta x} \frac{\sqrt{(x'-x_a)^2+\frac{1}{4}\Delta y^2+\frac{1}{4}\Delta z^2}}{(x'-x_a)^2+\frac{1}{4}\Delta z^2}\, \mathrm{d}x',\\
 l_{23}&=\Delta x \int_{z_k-\frac{1}{2}\Delta z}^{z_k+\frac{1}{2}\Delta z} \frac{\sqrt{\Delta x^2+\frac{1}{4}\Delta y^2+(z'-z_k)^2}}{\Delta x^2+(z'-z_k)^2}\, \mathrm{d}z',
\end{split}
\end{equation*}
and due to the symmetry of the integrand $ \bar {\bf r} \varphi(\bar r)$ on $xz$ plane
$$l_{24}=l_{22}.$$ 
So finally  we have $$s_2=\frac{1}{2}\Delta y (l_{21}+2 l_{22}+l_{23}).$$
Due to the symmetry of ${\bf r}$ in $V_{i,j,k}$ along $y$ direction, we have $$s_5=s_2.$$
The calculation of  $s_3$ is similar to the one of $s_2$ with the final result
$$s_3=\frac{1}{2}\Delta z (l_{31}+2 l_{32}+l_{33}),$$
where 
\begin{equation*}
\begin{split}
l_{31}&=-\frac{1}{2} \Delta z \pi,\\
 l_{32}&=\frac{1}{2}\Delta y \int_{x_a}^{x_a+\Delta x} \frac{\sqrt{(x'-x_a)^2+\frac{1}{4}\Delta z^2+\frac{1}{4}\Delta y^2}}{(x'-x_a)^2+\frac{1}{4}\Delta y^2}\, \mathrm{d}x',\\
 l_{33}&=\Delta x \int_{y_j-\frac{1}{2}\Delta y}^{y_j+\frac{1}{2}\Delta y} \frac{\sqrt{\Delta x^2+\frac{1}{4}\Delta z^2+(y'-y_j)^2}}{\Delta x^2+(y'-y_j)^2}\, \mathrm{d}y'.
\end{split}
\end{equation*}
Also due to the symmetry of ${\bf r}$ in $V_{i,j,k}$ along $z$ direction, we have $$s_6=s_3.$$
The only surface integral remaining to be calculated is the one over $S_4.$ On this surface we have $${\bf r}=(\Delta x, y'-y_j, z'-z_k),$$ and 
$${\bf n}_4=(1,0,0),$$
so that
\begin{equation*}
s_4=\Delta x\iint_{S_4} \frac{1}{\sqrt{\Delta x^2+(y'-y_j)^2+(z'-z_k)^2}} \, \mathrm{d} S.\label{v1s4}
\end{equation*}
Defining $$\bar {\bf r}=(y'-y_j,z'-z_k)$$ and $$\bar r=|\bar {\bf r}|,$$
we seek a function $\varphi(\bar r)$ that satisfies
$$\nabla \cdot (\bar {\bf r} \varphi(\bar r))=\frac{1}{\sqrt{\bar r^2+\Delta x^2}}.$$
This equation can be written in the form
$$2 \varphi(\bar r)+\bar r\varphi'(\bar r)=\frac{1}{\sqrt{\bar r^2+\Delta x^2}}.$$
Solving the above equation gives
\begin{equation*}\label{f1s4}
\varphi(\bar r)=\frac{\sqrt{\bar r^2+\Delta x^2}}{\bar r^2} .
\end{equation*}
\begin{figure}[h!]
\centering
\includegraphics[width=6cm]{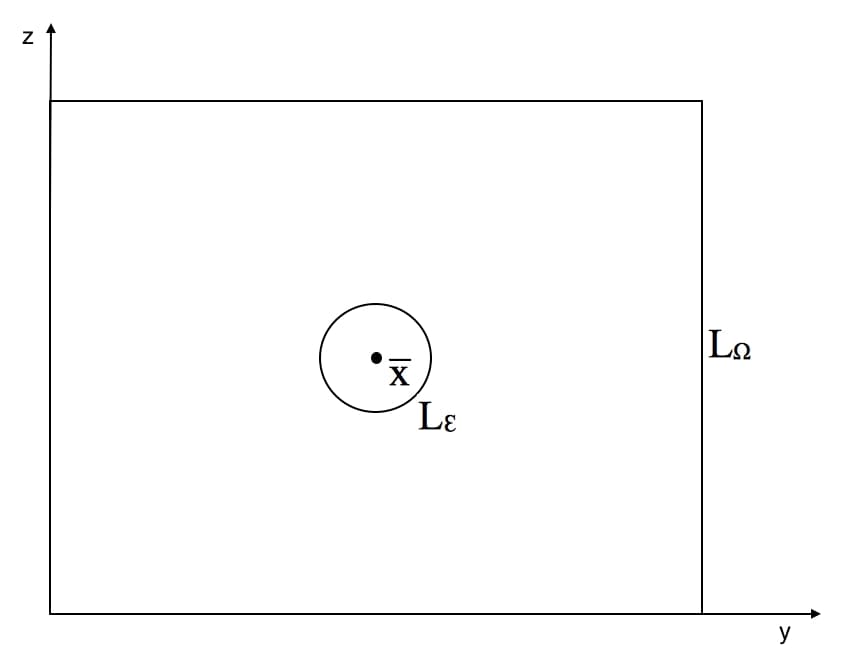}
\caption{Surface $S_4$ }
\label{s4}
\end{figure}
Applying the divergence theorem, we have 
\begin{equation*}
s_4=\Delta x\lim_{\epsilon \rightarrow 0}\int_{L_\epsilon}\frac{\sqrt{\bar r^2+\Delta x^2}}{\bar r^2}\bar {\bf r}\cdot\bar  {\bf n}_\epsilon\, \mathrm{d} L+\Delta x \int_{L_\Omega}\frac{\sqrt{\bar r^2+\Delta x^2}}{\bar r^2}\bar {\bf r}\cdot\bar  {\bf n}_\Omega\, \mathrm{d} L,
\end{equation*}
where $L_\epsilon $ is a circle with radius $\epsilon$ centered at point ${\bar{ \bf x}}=(y_j,z_k),$ and $L_\Omega $ is the four edges of surface $S_4$.
$\bar  {\bf n}_\epsilon$ is the unit normal vector of  $L_\epsilon$ and $\bar  {\bf n}_\Omega$ is the unit normal vector of  $L_\Omega,$ as shown in figure \ref{s4}.

For the integral over $L_\epsilon,$ we write
$$\bar {\bf r}= \epsilon (\cos \theta , \sin \theta), $$ 
and $$\bar {\bf n}_\epsilon= - (\cos \theta , \sin \theta),$$
 then 
\begin{equation}\label{h1s41}
\begin{split}
&\lim_{\epsilon \rightarrow 0} \int_{L_\epsilon}\frac{\sqrt{\bar r^2+\Delta x^2}}{\bar r^2}\bar {\bf r}\cdot \bar {\bf n}_\epsilon \, \mathrm{d} L\\
&=-\lim_{\epsilon \rightarrow 0}\int_{0}^{2\pi} \epsilon (\cos \theta , \sin \theta)\cdot  (\cos \theta , \sin \theta ) \frac{\sqrt{\epsilon^2+\Delta x^2}}{\epsilon^2} \epsilon \, \mathrm{d} \theta\\
&=-2\Delta x \pi.
\end{split}
\end{equation}
For the integral over $L_\Omega,$ there is no singularity anymore and this leads to
\begin{equation}\label{h1s42}
\begin{split}
&\int_{l_\Omega}\frac{\sqrt{\bar r^2+\Delta x^2}}{\bar r^2}\bar {\bf r}\cdot \bar {\bf n}_\Omega \, \mathrm{d} L\\
&=2l_{41}+2l_{42},
\end{split}
\end{equation}
with 
\begin{equation*}
\begin{split}
l_{41}= \frac{1}{2}\Delta y \int_{z_k-\Delta z}^{z_k+\Delta z}\frac{\sqrt{\Delta x^2+\frac{1}{4}\Delta y^2+(z'-z_k)^2}}{\frac{1}{4}\Delta y^2+(z'-z_k)^2}\, \mathrm{d}z',
\end{split}
\end{equation*}
and 
\begin{equation*}
\begin{split}
l_{42}=\Delta z \int_{y_j-\Delta y}^{y_j+\Delta y}\frac{\sqrt{\Delta x^2+\frac{1}{4}\Delta z^2+(y'-y_j)^2}}{\frac{1}{4}\Delta z^2+(y'-y_j)^2}\, \mathrm{d}y'.
\end{split}
\end{equation*}
Summing up (\ref{h1s41}) and (\ref{h1s42}), we obtain,
$$s_4=2\Delta x(l_{41}+l_{42}- \pi \Delta x).$$
We then finally get the following expression for $f_1$
\begin{equation*}
\begin{split}
f_1&=\frac{\Delta y}{2}(l_{21}+2l_{22}+l_{23})+\frac{\Delta z}{2}(l_{31}+2l_{32}+l_{33})\\
&+\Delta x (l_{41}+l_{42}- \pi \Delta x).
\end{split}
\end{equation*} 
All the line integrals $l_{21}$ etc are non-singular and can be calculated accurately using numerical integration.


\section{Summary}\label{summary3d}
In this paper we have, by considering 3D light scattering, discussed some important issues that we believe will be generic for numerical implementations of the EOS formulation for wave scattering. We have shown that the numerical instabilities can be thought as arising separately from the domain part and the boundary update part of the algorithm. We have argued that the instability arising from the boundary part of the algorithm is strongly related to the late time instability noted earlier while solving antenna problems using TBEM. We find that our version of the late time instability can be completely removed by suitably chosen material values, in particular the jump in material values at the boundary of the scattering object should not be too severe. In the limit where the material parameters simulate the properties of highly conductive metallic surfaces, we observe that our version of the late time instability is always present. Thus the instability interval vanishes in this limit. We take this as an indicator that for situations like in antenna theory, the late time instability should always be present, which it is. We are now aware of work where it has been noted that the instability can be removed by manipulating the material parameters defining the scattering objects. The EOS formulation gives thus different window into the late time instability that might be useful.

We have in our discretization used explicit methods. It would not be easy, but we  believe that it is possible to do a fully implicit method for  the EOS formulation, such an approach might remove all instabilities, which is the ultimate goal both for TBEM and for our EOS formulation.

In this paper we have also discussed how to calculate  singular volume and surface integrals for light scattering. The reason for including this discussion is that we think the type of singular integrals we discuss are generic for the singular integrals that will arise while calculating wave scattering using the EOS approach. 
\begin{appendices}
\section{Matrix elements}\label{Amatrix}
In this section we detail the entries of the updating matrix $M$ in (\ref{3dmatrixequation}) where $Q$ is a vector containing the components of the electric field and the magnetic field at all points of the grid with a size $6\times N_x \times N_y \times N_z$, where $N_x, $ $N_y$ and $N_z$ are the number of grid points in the  $x,$ $ y$ 
and $z$ directions.
To simplify the writing, we  denote 
\begin{equation*}
\Lambda_1=N_x\times N_y\times N_z,
\end{equation*}
\begin{equation*}
\Lambda_2= N_y\times N_z,
\end{equation*}
\begin{equation*}
\Lambda_3= 6\Lambda_1,
\end{equation*}
\begin{equation*}
\Gamma_1=\Lambda_2 i+N_z j+k,
\end{equation*}
\begin{equation*}
\Gamma_2=\Lambda_1+\Gamma_1,
\end{equation*}
\begin{equation*}
\Gamma_3=2\Lambda_1+\Gamma_1,
\end{equation*}
\begin{equation*}
\Gamma_4=3\Lambda_1+\Gamma_1, 
\end{equation*}
\begin{equation*}
\Gamma_5=4\Lambda_1+\Gamma_1, 
\end{equation*}
\begin{equation*}
\Gamma_6=5\Lambda_1+\Gamma_1.
\end{equation*}
Thus $Q$ is expressed by
\begin{equation*}
Q=
\left(\!
    \begin{array}{c}
   {[e_{1,i,j,k}]}_{\Lambda_1}\\{[e_{2,i,j,k}]}_{\Lambda_1}\\{[e_{3,i,j,k}]}_{\Lambda_1}\\{[b_{1,i,j,k}]}_{\Lambda_1}\\{[b_{2,i,j,k}]}_{\Lambda_1}\\{[b_{3,i,j,k}]}_{\Lambda_1}
    \end{array}
  \!\right)^{n+1}=\left(\!
    \begin{array}{c}
   {[Q_{\Gamma_1}]}_{\Lambda_1}\\{[Q_{\Gamma_2}]}_{\Lambda_1}\\{[Q_{\Gamma_3}]}_{\Lambda_1}\\{[Q_{\Gamma_4}]}_{\Lambda_1}\\{[Q_{\Gamma_5}]}_{\Lambda_1}\\{[Q_{\Gamma_6}]}_{\Lambda_1}
    \end{array}
  \!\right)^{n+1},
\end{equation*}
where  ${[e_{1,i,j,k}]}_{\Lambda_1}$ represents the vector containing the components of the electric field $e_1$ at all points of the grid indexing in $k, j, i$ order. ${[e_{2,i,j,k}]}_{\Lambda_1}$ and so on follow the same rule. 
Due to the complexity of the matrix, here we only illustrate the entries of the rows of $M$ corresponding to the components $Q_{\Gamma_1}^{n+1}.$ Other entries of the matrix can be expressed in the same way. 

After applying the Lax-Wendroff method, we have
\begin{equation*}
\begin{split}
e_{1,i,j,k}^{n+1}&=e_{1,i,j,k}^{n}+w_1 (e_{1,i,j,k}^n)_{yy}+w_1 (e_{1,i,j,k}^n)_{zz}-w_1 (e_{2,i,j,k}^n)_{xy}\\
&-w_1 (e_{3,i,j,k}^n)_{xz}+w_2(b_{3,i,j,k}^n)_y-w_2(b_{2,i,j,k}^n)_z,
\end{split}
\end{equation*}
$$\Downarrow$$
\begin{equation}\label{onerowmatrix}
\begin{split}
Q_{\Gamma_1}^{n+1}&=Q_{\Gamma_1}+w_1 (Q_{\Gamma_1})_{yy}+w_1 (Q_{\Gamma_1})_{zz}-w_1 (Q_{\Gamma_2})_{xy}\\
&-w_1 (Q_{\Gamma_3})_{xz}+w_2(Q_{\Gamma_6})_{y}-w_2(Q_{\Gamma_5})_{z},
\end{split}
\end{equation}
where 
\begin{equation*}
\begin{split}
w_1&=\frac{c^2\Delta t^2}{2},\\
w_2&=c^2\Delta t.
\end{split}
\end{equation*}
The coefficients of the right side of the equation (\ref{onerowmatrix}) are corresponding to the  $\Gamma_1$-th row of the matrix $M$ and the values of them are depended on the values of $i, j $ and $k.$ 
In order to have a compact and uniform expressions, we write 
\begin{equation*}
\begin{split}
&(Q_{\Gamma_6})_{y}=\frac{1}{ \Delta y}(\xi_{-2}Q_{\kappa_{-2}} +\xi_{-1} Q_{\kappa_{-1}}+\xi Q_{\kappa}+\xi_{1}Q_{\kappa_{1}} +\xi_{2} Q_{\kappa_{2}}),\\
&(Q_{\Gamma_1})_{yy}=\frac{1}{ (\Delta y)^2}(\delta_{-2}Q_{\chi_{-2}} +\delta_{-1} Q_{\chi_{-1}}+\delta Q_{\Gamma_1}+\delta_{1}Q_{\chi_{1}} +\delta_{2} Q_{\chi_{2}}),\\
&(Q_{\Gamma_2})_{xy}=\frac{1}{  3 \Delta x \Delta y}(\omega_{-4}Q_{\Upsilon_{-4}} +\omega_{-3} Q_{\Upsilon_{-3}}+\omega_{-2}Q_{\Upsilon_{-2}} +\omega_{-1} Q_{\Upsilon_{-1}}+\omega Q_{\Upsilon}\\
&\qquad \qquad +\omega_{1}Q_{\Upsilon_{1}} +\omega_{2} Q_{\Upsilon_{2}}+\omega_{3}Q_{\Upsilon_{3}} +\omega_{4} Q_{\Upsilon_{4}}),
\end{split}
\end{equation*}
where 
\begin{align*}
    &\chi_{1}=\Gamma_1+N_z,  &  &\chi_{2}=\Gamma_1+2N_z,& &\chi_{-1}=\Gamma_1-N_z, \\
 & \chi_{-2}=\Gamma_1-2N_z, & & \kappa=\Gamma_6, & &\kappa_{1}=\Gamma_6+N_z,  \\
 &\kappa_{2}=\Gamma_6+2N_z, & &\kappa_{-1}=\Gamma_6-N_z,  & &\kappa_{-2}=\Gamma_6-2N_z,\\
& \eta=N_y, & &\Upsilon=\Gamma_{2},& &\Upsilon_{-4}=\Gamma_{2}-\Lambda_2-N_z,\\
&\Upsilon_{-3}=\Gamma_{2}-\Lambda_2,& &\Upsilon_{-2}=\Gamma_{2}-\Lambda_2+N_z, & &\Upsilon_{-1}=\Gamma_{2}-N_z,\\
&\Upsilon_{1}=\Gamma_{2}+N_z,& &\Upsilon_{2}=\Gamma_{2}+\Lambda_2-N_z,& &\Upsilon_{3}=\Gamma_{2}+\Lambda_2,\\
& \Upsilon_{4}=\Gamma_{2}+\Lambda_2+N_z.
\end{align*}
The expressions for $(Q_{\Gamma_5})_{z},$ $(Q_{\Gamma_1})_{zz}$ and $(Q_{\Gamma_3})_{xz}$ have the same forms as $(Q_{\Gamma_6})_{y},$ $(Q_{\Gamma_1})_{yy}$ and $(Q_{\Gamma_2})_{xy}$ respectively, but with
\begin{align*}
    &\chi_{1}=\Gamma_1+1, &  &\chi_{2}=\Gamma_1+2,& &\chi_{-1}=\Gamma_1-1, \\
 & \chi_{-2}=\Gamma_1-2, & &\kappa=\Gamma_5, & &\kappa_{1}=\Gamma_5+1,  \\
 &\kappa_{2}=\Gamma_5+2,& &\kappa_{-1}=\Gamma_5-1,  & &\kappa_{-2}=\Gamma_5-2,\\
& \eta=N_z, & &\Upsilon=\Gamma_{3},& &\Upsilon_{-4}=\Gamma_{3}-\Lambda_2-1,\\
&\Upsilon_{-3}=\Gamma_{3}-\Lambda_2,& &\Upsilon_{-2}=\Gamma_{3}-\Lambda_2+1, & &\Upsilon_{-1}=\Gamma_{3}-1,\\
&\Upsilon_{1}=\Gamma_{3}+1,& &\Upsilon_{2}=\Gamma_{3}+\Lambda_2-1,& &\Upsilon_{3}=\Gamma_{3}+\Lambda_2,\\
&\Upsilon_{4}=\Gamma_{3}+\Lambda_2+1.
\end{align*}
After discussing the locations of $i, j$ and $k,$ the values of the coefficients are listed in table \ref{table1} and table \ref{table2}.
\begin{center}
\captionof{table}{($\frac{\partial}{\partial y}$,$\frac{\partial^2}{\partial y^2}$) or ($\frac{\partial}{\partial z}$,$\frac{\partial^2}{\partial z^2}$) related coefficients}
\label{table1}
    \begin{tabular}{| l | l | l | l | l | l | l | l | l | l |l|}
    \hline
    j or k &$\delta_{-2}$ & $\delta_{-1}$ & $\delta$ &$\delta_1$ & $\delta_2$&$\xi_{-2}$ &$\xi_{-1} $& $\xi$ & $\xi_1$ &$\xi_2$ \\ \hline
    0 & 0& 0&-5 & 2  & -1/5&0 &0&1/2 & 2/3 &-1/10 \\ \hline
    $\eta$-1 &-1/5 &2 & -5 & 0&0 & 1/10&-2/3&-1/2 & 0 &0 \\ \hline
   [1,$\eta$-2] &  0  &1 & -2 &1&0  &0&-1/2& 0 & 1/2 &0 \\ 
    \hline
    \end{tabular}
\end{center}
\begin{table*}
 \caption{$\frac{\partial^2}{\partial x \partial y}$ or $\frac{\partial^2}{\partial x \partial z}$ related coefficients}
\label{table2}
\begin{tabularx}{\textwidth}{l|l|lllllllll}
\toprule
i      & j or k& $\omega_{-4}$ & $\omega_{-3}$ & $\omega_{-2}$ & $\omega_{-1}$ & $\omega$ & $\omega_1$ & $\omega_2$ & $\omega_3$& $\omega_4$ \\ \hline
\multirow{4}{*}{0}         &0    & 0  &  0  &0 &0  & 9    & -5     &0   & -5    &1 \\ 
        & $\eta$-1   &  0&   0  &  0 &5    &-9     &  0  & -1    & 5  &0\\
        & [1,$\frac{\eta}{2}$)   & 0  &   0  &  0 & 0   & 3     & -3   & 1    & -1  &0\\ 
        & $[\frac{\eta}{2},\eta-2]$   &  0 &0 &0&3    & -3   &  0   &  0   & 1   & -1   \\\hline 
\multirow{4}{*}{  $N_x$-1}         &0 &  0  &  5 & -1   & 0 &- 9    & 5     & 0  & 0   &0 \\ 
        & $\eta$-1   &1  &  -5   & 0  &-5    &9     &0    &  0   & 0  &0\\
        & [1,$\frac{\eta}{2}$)   &  -1 & 1    &0   & 0   &- 3     & 3   &   0  & 0 &0\\ 
        & $[\frac{\eta}{2},\eta-2]$   & 0  & -1&1&-3    & 3   & 0    &  0   &  0  & 0   \\\hline
\multirow{3}{*}{[1,$N_x$-3]}         &0    &  0 & 0   &1 &0  & 3    & -1     & 0  & -3    &0 \\ 
        & $\eta$-1   & -1 &   0  &  0 &1    &-3     & 0   &    0 & 3  &0\\ 
        & $[1,\eta-2]$   &  -3/4 &0 &3/4&  0  &  0 &   0  &3/4     & 0   & -3/4   \\\hline
\multirow{3}{*}{$N_x$-2}         &0    & 0  &  3  &0 & 0 &-3    &1     &0   & 0    &-1 \\ 
        & $\eta$-1   & 0 &  -3   &  0 &-1    & 3    & 0   &  1   & 0  &0\\
        & [1,$\eta$-2] &  -3/4 &0 &3/4&  0  & 0  & 0    &3/4     & 0   & -3/4   \\
\bottomrule
\end{tabularx}
\end{table*}

For example, if $i=0, $ $j=0$ and $k=0$, the entries of the $\Gamma_1$-th row of the matrix $M$ are the following
\begin{align*}
    &M_{\Gamma_1,\Gamma_1}=1-5u_1-5v_1,  &  &M_{\Gamma_1,\Gamma_1+N_z}=2u_1, & &M_{\Gamma_1,\Gamma_1+2 N_z}=-\frac{1}{5}u_1, \\
 & M_{\Gamma_1,\Gamma_6}=\frac{1}{2}u_2,  & & M_{\Gamma_1,\Gamma_6+N_z}=\frac{2}{3}u_2, & &M_{\Gamma_1,\Gamma_6+2 N_z}=-\frac{1}{10}u_2,  \\
 &M_{\Gamma_1,\Gamma_2}=9u_3, & &M_{\Gamma_1,\Gamma_2+N_z}=-5u_3,  & &M_{\Gamma_1,\Gamma_2+\Lambda_2}=-5u_3, \\
& M_{\Gamma_1,\Gamma_2+\Lambda_2+N_z}=u_3, & &M_{\Gamma_1,\Gamma_1+1}=2v_1, & &M_{\Gamma_1,\Gamma_1+2 }=-\frac{1}{5}v_1,\\
&M_{\Gamma_1,\Gamma_5}=\frac{1}{2}v_2,& & M_{\Gamma_1,\Gamma_5+1}=\frac{2}{3}v_2,& &M_{\Gamma_1,\Gamma_5+2 }=-\frac{1}{10}v_2,\\
&M_{\Gamma_1,\Gamma_3}=9v_3,& &M_{\Gamma_1,\Gamma_3+1}=-5v_3,& &M_{\Gamma_1,\Gamma_3+\Lambda_2}=-5v_3,\\
&M_{\Gamma_1,\Gamma_3+\Lambda_2+1}=v_3,
\end{align*}
otherwise $M_{\Gamma_1, *}=0$  and where 
\begin{align*}
    &u_1=\frac{w_1}{(\Delta y)^2}, &  &u_2=\frac{w_2}{\Delta y},& &u_3=\frac{w_1}{3 \Delta x\Delta y}, \\
 & v_1=\frac{w_1}{(\Delta z)^2},& &v_2=\frac{w_2}{\Delta z},& &v_3=\frac{w_1}{3 \Delta x\Delta z}.  
\end{align*}

\section{Singular integrals}\label{Singularity}
In this section, we detail the calculations of other types of singular integrals involved in the EOS formulations of 3D Maxwell's equations, denoted by ${\bf f}_2,$ $ {\bf f}_3,$ $ g_1,$ ${\bf g}_2,$ ${\bf g}_3$ in \cite{Aihua2}. 
The techniques are similar with the calculating of $f_1$ in section \ref{Asingularity}. The geometry  is llustrated in figure \ref{box}.  
\subsection{Calculation of ${\bf f}_2$}
\begin{equation}\label{f2}
 {\bf f}_2=\iiint_{V_{i,j,k}} \frac{{\bf x}'-{\bf x}_p}{|{\bf x}'-{\bf x}_p|^2}\,\mathrm{d}V.
\end{equation}
The components of the integration variable in (\ref{f2}) are given by
$$
{\bf x'}=(x', y', z'),
$$
and let us introduce the quantity
$$
{\bf r}={\bf x}'-{\bf x}_p,
$$
with $r=|{\bf r}|.$

We want to apply the divergence theorem on (\ref{f2}), and therefore need to find a function $\varphi(r)$ that satisfies 
 $$\nabla \cdot ({\bf r} {\bf r} \varphi(r))=\frac{{\bf r}}{r^2},$$
or equivalently
$$\nabla \cdot ({\bf r} {\bf r}) \varphi(r)+{\bf r} {\bf r}\cdot \nabla\varphi(r)=\frac{{\bf r} }{r^2} .$$
Solving the above equation, we get
$$ \varphi(r)=\frac{1}{2r^2}.$$
Because of the singularity on surface $S_1,$  we can not apply the divergence theorem directly. However we can write 
\begin{equation*}\label{f2i}
\begin{split}
{\bf f}_2= \frac{1}{2 }(\sum_{m=2}^6 \iint_{S_m} \frac{{\bf r}{\bf r}}{r^2}\cdot {\bf n}_m\, \mathrm{d} S+\lim_{\epsilon \rightarrow 0}\iint_{S_\epsilon}\frac{{\bf r}{\bf r}}{ r^2}\cdot {\bf n}_\epsilon\, \mathrm{d}S+\iint_{S_\Omega}\frac{{\bf r}{\bf r}}{r^2}\cdot {\bf n}_1\,\mathrm{d}S),
\end{split}
\end{equation*}
where $S_\epsilon$ is a hemispherical surface of radius $\epsilon$ centered at ${\bf x}_p$ and  $S_\Omega$ is the rest of the surface $S_1$ with a disk of radius $\epsilon$ around ${\bf x}_p$ has been removed. ${\bf n}_\epsilon$ is the unit normal vector on $S_\epsilon,$ pointing out of $V_{i,j,k}.$ 
${\bf n}_m$ is the unit normal vector on $S_m,$ pointing out of $V_{i,j,k}.$ 

For the integral over $S_\Omega$, we have
$${\bf r}=(0, y'-y_j, z'-z_k),$$ 
and
$${\bf n}_1=(-1,0,0),$$ 
thus we get
\begin{equation*}\label{h2s11}
\begin{split}
\iint_{S_\Omega}\frac{{\bf r}{\bf r}}{ r^2}\cdot {\bf n}\,\mathrm{d}S=(0, 0, 0). 
\end{split}
\end{equation*}
For the integral over surface $S_\epsilon$, we use the spherical coordinate system, 
$${\bf r}=\epsilon (\cos \theta \sin \varphi, \sin \theta \sin \varphi, \cos \varphi),$$
and 
$${\bf n}_\epsilon= (\cos \theta \sin \varphi, \sin \theta \sin \varphi, \cos \varphi),$$
where $\epsilon, \varphi, \theta$ are  respectively the radial distance, polar angle and azimuthal angle,
so that
\begin{equation*}\label{h2s12}
\begin{split}
\iint_{S_\epsilon}\frac{1}{ r^2} {\bf r}{\bf r}\cdot {\bf n}_\epsilon\,\mathrm{d}S=&\lim_{\epsilon \rightarrow 0}\frac{1}{\epsilon^2}\int_{0}^{2\pi}\int_{-\frac{\pi}{2}}^{\frac{\pi}{2}} \epsilon (\cos \theta \sin \varphi, \sin \theta \sin \varphi, \cos \varphi)\\
&\epsilon (\cos \theta \sin \varphi, \sin \theta \sin \varphi, \cos \varphi) \\
&\cdot  (\cos \theta \sin \varphi, \sin \theta \sin \varphi, \cos \varphi)\\
& \epsilon^2 \sin \varphi \, \mathrm{d} \theta \, \mathrm{d} \varphi=(0,0,0).
\end{split}
\end{equation*}
Defining 
\begin{equation*}
\begin{split}
s_m=  \iint_{S_m} \frac{{\bf r}{\bf r}}{r^2}\cdot {\bf n}_m\, \mathrm{d} S,
\end{split}
\end{equation*}
$f_2$ can be written as
\begin{equation*}
\begin{split}
{\bf f}_2= \frac{1}{2 } \sum_{m=2}^6 s_m.
\end{split}
\end{equation*}
Due to the symmetry of ${\bf r}$ along $y$ and $z$ directions in $V_{i,j,k}$, we have $$s_5= s_2$$ and $$s_6=s_3.$$
Thus $f_2$ can be written as,
\begin{equation*}
\begin{split}
{\bf f}_2&=\frac{1}{2}(\iint_{S_2}\Delta y \frac{(x'-x_a,0,z'-z_k)}{(x'-x_a)^2+(z'-z_k)^2+\frac{1}{4}\Delta y^2}\, \mathrm{d}x'\, \mathrm{d}z'\\
&+\iint_{S_3}\Delta z \frac{(x'-x_a,y'-y_j,0)}{(x'-x_a)^2+(y'-y_j)^2+\frac{1}{4}\Delta z^2}\, \mathrm{d}x'\, \mathrm{d}y'\\
&+\iint_{S_4}\Delta x^2 \frac{(1,0,0)}{\Delta x^2+(y'-y_j)^2+(z'-z_k)^2}\, \mathrm{d}y'\, \mathrm{d}z').
\end{split}
\end{equation*}
For computation simplicity, we define  
$$\bar{s}_2=\iint_{S_2} \frac{(x'-x_a,z'-z_k)}{(x'-x)^2+(z'-z_k)^2+\frac{1}{4}\Delta y^2}\, \mathrm{d}x'\, \mathrm{d}z',$$
$$\bar{s}_3=\iint_{S_3} \frac{(x'-x_a,y'-y_j)}{(x'-x_a)^2+(y'-y_j)^2+\frac{1}{4}\Delta z^2}\, \mathrm{d}x'\, \mathrm{d}y',$$
$$\bar{s}_4=\iint_{S_4} x^2 \frac{1}{\Delta x^2+(y'-y_j)^2+(z'-z_k)^2}\, \mathrm{d}y'\, \mathrm{d}z'.$$
Thus
for the calculations of $\bar{s}_2$ and $\bar{s}_3$, 
we consider a general form
\begin{equation}\label{f23}
\iint_{S}\frac{\bar {\bf r}}{\bar r^2+A^2}\, \mathrm{d}S,
\end{equation}
 where $A$ is a constant, $\bar {\bf r}$ is a 2-component vector on surface S, and $\bar r=|\bar {\bf r}|.$
 We want to apply the divergence theorem on (\ref{f23}), therefore we need to find a function $\varphi(\bar r)$ that satisfies 
$$\nabla \cdot (\bar{\bf r} \bar{\bf r} \varphi(\bar r))=\frac{\bar{\bf r}}{\bar r^2+A^2}.$$ 
Solving the above equation, we get
\begin{equation*}\label{f2s2}
\varphi(\bar r)=-\frac{A\tan^{-1}(\frac{\bar r}{A})}{\bar r^3}+\frac{1}{\bar r^2}.
\end{equation*}
For $\bar{s}_2,$ 
$$S=S_2, $$
$$A=\frac{1}{2}\Delta y,$$  and 
$$\bar{\bf r}=(x'-x_a,z'-z_k),$$ 
and because of the singularity on $S_2,$ we can not use the divergence theorem directly, however
we can write
\begin{equation*}
\begin{split}
 \bar{s}_2=\sum_{n=2}^4\int_{L_{2n}}\varphi(\bar r)\bar {\bf r}\bar {\bf r}\cdot\bar  {\bf n}_n\, \mathrm{d}L+\lim_{\epsilon\rightarrow 0}\int_{L_\epsilon} \varphi(\bar r)\bar {\bf r}\bar  {\bf r}\cdot \bar  {\bf n}_\epsilon\, \mathrm{d}L+\int_{L_\Omega} \varphi(\bar r)\bar {\bf r}\bar  {\bf r}\cdot\bar  {\bf n}_1\, \mathrm{d}L,
\end{split}
\end{equation*}
where $L_{2n}$ are edges of $S_2.$ 
 $L_\epsilon$ is a semicircle with radius $\epsilon$ centered at point ${\bar { \bf x}}$ and $L_\Omega$ is the rest of $L_{21}$. 
$\bar {\bf n}_\epsilon$ is the unit normal of $L_{\epsilon},$ pointing out of $S_2.$ 
$\bar {\bf n}_n$ is the unit normal of $L_{2n},$ pointing out of $S_2$.
Geometry is illustrated in figure \ref{s2}.

For the integral over $L_\Omega,$ 
we have $$\bar {\bf r}=(0, z'-z_k)$$ and $$\bar {\bf n}_1=(-1,0),$$ so that 
$$\int_{L_\Omega}\varphi(\bar r)\bar {\bf r}\bar  {\bf r}\cdot\bar  {\bf n}_1\, \mathrm{d}L=(0,0).$$
For the integral over $L_\epsilon,$ using the polar coordinates, we have $${\bf r}=\epsilon (\cos \theta , \sin \theta),$$
and $${\bf n}_\epsilon=-(\cos \theta , \sin \theta),$$
so that
\begin{equation*}
\begin{split}
\int_{L_\epsilon}\varphi(\bar r)\bar {\bf r}\bar  {\bf r}\cdot\bar  {\bf n}_\epsilon\, \mathrm{d}L&=-\lim_{\epsilon\rightarrow 0}\int_{-\frac{\pi}{2}}^{\frac{\pi}{2}} -\epsilon^3(\cos\theta,\sin \theta)(\cos\theta,\sin \theta)\\
&(-\frac{\frac{1}{2}\Delta y\tan^{-1}(\frac{\epsilon}{\frac{1}{2}\Delta y})}{\epsilon^3 }+\frac{1}{\epsilon^2}) \cdot (\cos\theta,\sin \theta)\, \mathrm{d}\theta\\
&=(0,0).
\end{split}
\end{equation*}
There is no singularity on $L_{22}, L_{23}$ and $L_{24},$
finally,
\begin{equation*}
\bar{s}_2=(I_1,0),
\end{equation*}
where 
\begin{equation*}
I_1=\int_{x_a}^{x_a+\Delta x}\Delta z (x'-x_a)\varphi(r_1)\, \mathrm{d}x'+\int_{z_k-\Delta z}^{z_k+\Delta z}\Delta x^2\varphi(r_2)\, \mathrm{d}z',
\end{equation*}
 with
$$ r_1=\sqrt{(x'-x_a)^2+\frac{1}{4}\Delta z^2},$$ and $$ r_2=\sqrt{\Delta x^2+(z'-z_k)^2}.$$ 
The calculation of  $\bar{s}_3$ is similar to the one of $\bar{s}_2$ with the final result
\begin{equation*}
\bar{s}_3=(I_2,0),
\end{equation*}
where
\begin{equation*}
I_2=\int_{x_a}^{x_a+\Delta x}\Delta y (x'-x_a)\varphi(r_3)\, \mathrm{d}x'+\int_{y_j-\Delta y}^{y_j+\Delta y}\Delta x^2\varphi(r_4)\, \mathrm{d}y',
\end{equation*}
with 
 $$ r_3=\sqrt{(x'-x_a)^2+\frac{1}{4}\Delta y^2},$$
and 
$$ r_4=\sqrt{\Delta x^2+(y'-y_j)^2}.$$

For the integral 
$\bar{s}_4,$  defining
\begin{equation*}
\bar{ \bf r}=(y'-y_j,z'-z_k),
 \end{equation*}
and $${\bar r}=|\bar{ \bf r}|,$$
we need to find a function $\varphi(\bar r)$ that satisfies
$$\nabla \cdot (\bar {\bf r}\varphi(\bar r))=\frac{1}{\bar r^2+\Delta x^2}.$$
Solving the above equation gives
\begin{equation*}\label{f2s4}
\varphi(\bar r)=\frac{\ln(\bar r^2+\Delta x^2)}{2\bar r^2}.
\end{equation*}
Because of the singularity at point ${\bar {\bf x}}=(y_j, z_k),$ we write
$$ \bar{s}_4=\lim_{\epsilon\rightarrow 0}\int_{L_{\epsilon}}\varphi(\bar r)\bar {\bf r}\cdot\bar  {\bf n}_\epsilon\,\mathrm{d}L+\int_{L_{\Omega}}\varphi(\bar r)\bar {\bf r}\cdot \bar {\bf n}_\Omega\,\mathrm{d}L,$$ 
where $L_{\epsilon}$ is a circle with radius $\epsilon$ centered at ${\bar {\bf x}}$ and $L_{\Omega}$ is the four edges of surface $S_4.$ $\bar  {\bf n}_\epsilon$ is the unit normal vector of $L_\epsilon,$ $ \bar {\bf n}_\Omega$ is the unit normal vector of $L_\Omega,$ as shown in figure \ref{s4}.

For the integral over $L_{\epsilon}$, we use the polar coordinates,
$$\bar {\bf r}=\epsilon (\cos\theta,\sin\theta),$$ 
and $${\bf n}_\epsilon=- (\cos\theta,\sin\theta),$$ 
then
\begin{equation*}
\begin{split}
\lim_{\epsilon\rightarrow 0}\int_{L_{\epsilon}}\varphi(\bar r)\bar {\bf r}\cdot\bar  {\bf n}_\epsilon\,\mathrm{d}L&=\lim_{\epsilon\rightarrow 0}\int_{0}^{2\pi}-\epsilon^2(\cos\theta,\sin\theta)\frac{\ln(\epsilon^2+\Delta x^2)}{2\epsilon^2}(\cos\theta,\sin\theta)\, \mathrm{d}\theta\\
&=-\pi\ln(\Delta x^2).
\end{split}
\end{equation*}
For the integral over $L_{\Omega},$ there is no singularity any more and 
this leads to
\begin{equation*}
\begin{split}
\int_{l_{\Omega}}\varphi(\bar r)\bar {\bf r}\cdot \bar {\bf n}_\Omega\,\mathrm{d}L&=\frac{1}{2}\Delta y\int_{z_k-\Delta z}^{z_k+\Delta z}\frac{\ln(\frac{1}{4}\Delta y^2+\Delta x^2+(z'-z_k)^2)}{\frac{1}{4}\Delta y^2+(z'-z_k)^2}\, \mathrm{d}z'\\
&+\frac{1}{2}\Delta z\int_{y_j-\Delta y}^{y_j+\Delta y}\frac{\ln(\frac{1}{4}\Delta z^2+\Delta x^2+(y'-y_j)^2)}{\frac{1}{4}\Delta z^2+(y'-y_j)^2}\, \mathrm{d}y'\\
&=I_3.
\end{split}
\end{equation*}
Combining all the above calculations, we finally get, 
$${\bf f}_2=\frac{1}{2}(\Delta y I_1+\Delta z I_2+\Delta x^2 (I_3-\pi\ln(\Delta x^2)), 0, 0).$$
All the line integrals $I_1,$ $I_2,$ $I_3$ are non-singular and can be calculated using numerical integration.\\

\subsection{Calculation of ${\bf f}_3$}
 \begin{equation}\label{f3}
{\bf f}_3=\iiint_{V_{i,j,k}}\frac{{\bf x}'-{\bf x}_p}{|{\bf x}'-{\bf x}_p|^3}\,\mathrm{d}V.
\end{equation}
The components of the integration variable in (\ref{f3}) are given by
$$
{\bf x'}=(x', y', z'),
$$
and let us introduce the quantity
$$
{\bf r}={\bf x}'-{\bf x}_p,
$$
with $r=|{\bf r}|.$

We want to apply the divergence theorem on (\ref{f3}), and therefore need to find a function $\varphi(r)$ that satisfies
$$\nabla \cdot ({\bf r} {\bf r} \varphi(r))=\frac{{\bf r}}{r^3},$$
thus we have
$$\nabla \cdot ({\bf r} {\bf r}) \varphi(r)+{\bf r} {\bf r}\cdot \nabla\varphi(r)=4{\bf r} \varphi(r)+r{\bf r} \varphi'(r)=\frac{{\bf r} }{r^3} .$$
Solving the above equation gives
\begin{equation*}
\varphi(r)=\frac{1}{r^3}.
\end{equation*}
Because of the singularity on surface $S_1$, we can not apply the divergence theorem directly, however we can write 
\begin{equation*}\label{f3i}
\begin{split}
{\bf f}_3&=\sum_{m=2}^6 \iint_{S_m} \frac{{\bf r}{\bf r}}{ r^3}\cdot {\bf n}_m\, \mathrm{d} S+\lim_{\epsilon \rightarrow 0}\iint_{S_\epsilon}\frac{{\bf r}{\bf r}}{ r^3}\cdot {\bf n}_\epsilon\,\mathrm{d}S+\iint_{S_\Omega}\frac{{\bf r}{\bf r}}{ r^3}\cdot {\bf n}_1\,\mathrm{d}S,
\end{split}
\end{equation*}
where $S_\epsilon$ is a hemispherical surface of radius $\epsilon$ centered at ${\bf x}_p$ and  $S_\Omega$ is the rest of the surface $S_1$ with a disk of radius $\epsilon$ around ${\bf x}_p$ has been removed. ${\bf n}_\epsilon$ is the unit normal vector on $S_\epsilon,$ pointing out of $V_{i,j,k}.$ 
${\bf n}_m$ is the unit normal vector on $S_m,$ pointing out of $V_{i,j,k}.$ 

For the integral over $S_\Omega$, we have
 $${\bf r}=(0, y'-y_j, z'-z_k),$$
and
$${\bf n}_1=(-1,0,0) ,$$
thus we get 
\begin{equation*}\label{f3s11}
\begin{split}
\iint_{S_\Omega}\frac{{\bf r}{\bf r}}{r^3}\cdot {\bf n}\,\mathrm{d}S=(0, 0, 0). 
\end{split}
\end{equation*}
For the integral over surface $S_\epsilon$, we use the spherical coordinate system, 
$${\bf r}=\epsilon (\cos \theta \sin \varphi, \sin \theta \sin \varphi, \cos \varphi),$$
and 
$${\bf n}_\epsilon= (\cos \theta \sin \varphi, \sin \theta \sin \varphi, \cos \varphi),$$
where $\epsilon, \varphi, \theta$ are  respectively the radial distance, polar angle and azimuthal angle,
so that
\begin{equation*}\label{f3s12}
\begin{split}
\iint_{S_\epsilon}\frac{1}{2 r^2} {\bf r}{\bf r}\cdot {\bf n}_\epsilon\,\mathrm{d}S=&\lim_{\epsilon \rightarrow 0}\int_{0}^{2\pi}\int_{-\frac{\pi}{2}}^{\frac{\pi}{2}}\frac{1}{2 r^2} {\bf r}{\bf r}\cdot {\bf n}_\epsilon \epsilon^2 \sin \varphi \,\mathrm{d} \theta \,\mathrm{d} \varphi=(0, 0, 0).
\end{split}
\end{equation*}
Defining 
\begin{equation*}
\begin{split}
s_m=\iint_{S_m} \frac{{\bf r}{\bf r}}{ r^3}\cdot {\bf n}_m\, \mathrm{d} S,
\end{split}
\end{equation*}
${\bf f}_3$ can be written as
\begin{equation*}
\begin{split}
{\bf f}_3=\sum_{m=2}^6 s_m.
\end{split}
\end{equation*}
Due to the symmetry of ${\bf r}$ in $V_{i,j,k}$ along $y $ and $z$ direction , we have  
$$s_2= s_5$$ and $$s_3=s_6.$$ So we have the following,
\begin{equation*}
\begin{split}
{\bf f}_3&=\iint_{S_2}\Delta y \frac{(x'-x_a,0,z'-z_k)}{((x'-x_a)^2+(z'-z_k)^2+\frac{1}{4}\Delta y^2)^{\frac{3}{2}}}\,\mathrm{d}x'\,\mathrm{d}z'\\
&+\iint_{S_3}\Delta z \frac{(x'-x_a,y'-y_j,0)}{((x'-x_a)^2+(y'-y_j)^2+\frac{1}{4}\Delta z^2)^{\frac{3}{2}}}\,\mathrm{d}x'\,\mathrm{d}y'\\
&+\iint_{S_4}\Delta x \frac{(\Delta x,0,0)}{(\Delta x^2+(y'-y_j)^2+(z'-z_k)^2)^{\frac{3}{2}}}\,\mathrm{d}y'\,\mathrm{d}z',
\end{split}
\end{equation*}
For computation simplicity, we define 
$${\bar s_2}=\iint_{S_2}\frac{(x'-x_a,z'-z_k)}{((x'-x_a)^2+(z'-z_k)^2+\frac{1}{4}\Delta y^2)^{3/2}}\,\mathrm{d}x'\,\mathrm{d}z',$$
$${\bar s_3}=\iint_{S_3} \frac{(x'-x_a,y'-y_j)}{((x'-x_a)^2+(y'-y_j)^2+\frac{1}{4}\Delta z^2)^{3/2}}\,\mathrm{d}x'\,\mathrm{d}y',$$
and 
$${\bar s_4}=\iint_{S_4} \frac{\Delta x}{(\Delta x^2+(y'-y_j)^2+(z'-z_k)^2)^{3/2}}\,\mathrm{d}y'\,\mathrm{d}z'.$$
Thus for the calculations of ${\bar s}_2$ and ${\bar s}_3$, we consider a general  form  
\begin{equation}\label{f33}
\iint_{S}\frac{\bar {\bf r}}{(\bar r^2+A^2)^{3/2}}\,\mathrm{d}S,
\end{equation}
where $\bar {\bf r}$ is a 2-component vector, $A$ is a constant and $\bar r=\bar {\bf r}.$ 
We want to apply the divergence theorem on (\ref{f33}), thus we need to find a function $\varphi(\bar r)$ that satisfies
$$\nabla \cdot (\bar {\bf r} \bar {\bf r} \varphi(\bar r))=\frac{\bar {\bf r}}{(\bar r^2+A^2)^{3/2}}.$$  
Solving the above equation,  we get
\begin{equation*}\label{f3s2}
\varphi(\bar r)=\frac{\log(\sqrt{\bar r^2+A^2}+\bar r)}{\bar r^3}-\frac{1}{\bar r^2\sqrt{\bar r^2+A^2}}.
\end{equation*}
For ${\bar s_2},$ 
$$S=S_2,$$
$$A=\frac{1}{2}\Delta y,$$ and $${\bar {\bf x}}=(x'-x_a,z'-z_k),$$
because of the singularity on $ S_2,$ we write
\begin{equation*}
\begin{split}
{\bar s_2}= \sum_{n=2}^4  \int_{L_{2n}}\varphi(\bar r){\bf r}\bar {\bf r}\cdot \bar {\bf n}_n\,\mathrm{d}L+\lim_{\epsilon \rightarrow 0}\int_{L_\epsilon} \varphi(\bar r)\bar {\bf r}\bar  {\bf r}\cdot \bar {\bf n}_\epsilon\,\mathrm{d}L+\int_{L_\Omega} \varphi(\bar r)\bar {\bf r} \bar {\bf r}\cdot \bar {\bf n}_1\,\mathrm{d}L,
\end{split}
\end{equation*}
where $L_{2n}$ are edges of $S_2.$ 
 $L_\epsilon$ is a semicircle with radius $\epsilon$ centered at point ${\bar { \bf x}}$ and $L_\Omega$ is the rest of $L_{21}$. 
$\bar {\bf n}_\epsilon$ is the unit normal of $L_{\epsilon},$ pointing out of $S_2.$ 
$\bar {\bf n}_n$ is the unit normal of $L_{2n},$ pointing out of $S_2$.
Geometry is illustrated in figure \ref{s2}.

For the integral over $L_\Omega,$ we have
$$\bar {\bf r}=(0, z'-z_k)$$ and $$\bar {\bf n}=(-1,0),$$ 
so that $$\int_{L_\Omega} \varphi(\bar r)\bar {\bf r}\bar  {\bf r}\cdot \bar {\bf n}_1\,\mathrm{d}L=(0,0).$$
For the integral over  $L_\epsilon,$ using the polar coordinates, we have
 $$\bar {\bf r}=\epsilon (\cos \theta , \sin \theta),$$
and 
$$\bar {\bf n}_\epsilon=-(\cos \theta , \sin \theta),$$
then
\begin{equation*}
\begin{split}
&\int_{L_\epsilon}  \varphi(\bar r)\bar {\bf r} \bar {\bf r}\cdot \bar {\bf n}_1\,\mathrm{d}L\\
&=\lim_{\epsilon\rightarrow 0}\int_{-\frac{\pi}{2}}^{\frac{\pi}{2}}-\epsilon^2(\cos\theta,\sin \theta)^3 (-\frac{1}{\epsilon^2(\epsilon^2+\frac{1}{4}\Delta y^2)^{\frac{1}{2}}}+\frac{\log(\sqrt{\epsilon^2+\frac{1}{4}\Delta y^2}+\epsilon)}{\epsilon^3})\,\mathrm{d}\theta\\
&=(-2\log(\frac{1}{2}\Delta y), 0),
\end{split}
\end{equation*}
There is no singularity on $L_{22}, L_{23}$ and $L_{24}$ any more,
finally,
$$ {\bar s_2}=(I_1-2\log(\frac{1}{2}\Delta y), 0),$$ 
where 
\begin{equation*}
I_1=\int_{x_a}^{x_a+\Delta x}\Delta z (x'-x_a)\varphi(r_1)\, \mathrm{d}x'+\int_{z_k-\Delta z}^{z_k+\Delta z}\Delta x^2\varphi(r_2)\, \mathrm{d}z',
\end{equation*}
 with
$$ r_1=\sqrt{(x'-x_a)^2+\frac{1}{4}\Delta z^2},$$ and $$ r_2=\sqrt{\Delta x^2+(z'-z_k)^2}.$$ 
The calculation of  $\bar{s}_3$ is similar to the one of $\bar{s}_2$ with the final result
\begin{equation*}
\bar{s}_3=(I_2-2\log(\frac{1}{2}\Delta z),0),
\end{equation*}
where
\begin{equation*}
I_2=\int_{x_a}^{x_a+\Delta x}\Delta y (x'-x_a)\varphi(r_3)\, \mathrm{d}x'+\int_{y_j-\Delta y}^{y_j+\Delta y}\Delta x^2\varphi(r_4)\, \mathrm{d}y',
\end{equation*}
with 
 $$ r_3=\sqrt{(x'-x_a)^2+\frac{1}{4}\Delta y^2},$$
and 
$$ r_4=\sqrt{\Delta x^2+(y'-y_j)^2}.$$

For the integral 
${\bar s_4},$  defining 
$$\bar {\bf r}=(y'-y_j, z'-z_k)$$
and $$\bar r=|\bar {\bf r}|,$$
we seek a function that satisfies 
$$\nabla \cdot (\bar {\bf r}\varphi(\bar r))=\frac{1}{(\bar r^2+\Delta x^2)^{3/2}}.$$
Solving this equation gives
\begin{equation*}\label{f3s4}
  \varphi(\bar r)=-\frac{1}{\bar r^2\sqrt{\bar r^2+\Delta x^2}}.
\end{equation*}
Because of the singularity at point ${\bar {\bf x}}=(y_j, z_k),$ we write
$$ {\bar s_4}=\int_{L_{\epsilon}}\varphi(\bar r)\bar {\bf r}\cdot \bar {\bf n}_\epsilon\,\mathrm{d}L+\int_{L_{\Omega}}\varphi(\bar r)\bar {\bf r}\cdot \bar {\bf n}_\Omega\,\mathrm{d}L,$$ where $L_{\epsilon}$ is a circle with radius $\epsilon$ centered at ${\bar {\bf x}}$ and $L_{\Omega}$ is the four edges of surface $S_4.$ $\bar  {\bf n}_\epsilon$ is the unit normal vector of $L_\epsilon,$ $ \bar {\bf n}_\Omega$ is the unit normal vector of $L_\Omega,$ as shown in figure \ref{s4}.

For the integral over $L_{\epsilon}$, we use the polar coordinates,
$$\bar {\bf r}=\epsilon (\cos\theta,\sin\theta),$$ 
and $${\bf n}_\epsilon=- (\cos\theta,\sin\theta),$$ 
then
\begin{equation*}
\lim_{\epsilon\rightarrow 0}\int_{L_{\epsilon}}\varphi(\bar r)\bar {\bf r}\cdot \bar {\bf n}_\epsilon\,\mathrm{d}L=\frac{2\pi}{\Delta x}.
\end{equation*}
For the integral on $L_{\Omega},$ there is no singularity any more
and this gives
\begin{equation*}
\begin{split}
&\int_{L_{\Omega}}\varphi(\bar r)\bar {\bf r}\cdot \bar {\bf n}_\Omega\,\mathrm{d}L\\
&=-\Delta y\int_{z_k-\Delta z}^{z_k+\Delta z}\frac{1}{(\frac{1}{4}\Delta y^2+(z'-z_k)^2)\sqrt{\Delta x^2+\frac{1}{4}\Delta y^2+(z'-z_k)^2}}\,\mathrm{d}z'\\
&-\Delta z\int_{y_j-\Delta y}^{y_j+\Delta y}\frac{1}{(\frac{1}{4}\Delta z^2+(y'-y_j)^2)\sqrt{\Delta x^2+\frac{1}{4}\Delta z^2+(y'-y_j)^2}}\,\mathrm{d}y'\\
&=I_3.
\end{split}
\end{equation*} 
Combining all the calculations above, we finally get 
\begin{equation*}
\begin{split}
{\bf f}_3&=(\Delta y I_1+\Delta z I_2-2 \Delta y \log(\frac{1}{2}\Delta y)-2 \Delta z \log(\frac{1}{2}\Delta z)\\
&+\Delta x^2 I_3+2\pi \Delta x , 0, 0).
\end{split}
\end{equation*}
All the line integrals $I_1,$ $I_2,$ $I_3$ are non-singular and can be calculated using numerical integration.
\subsection{Calculations of singular surface integrals}
When the observing point ${\bf x}_p$ and the integrating point are both located on the same integral surface, for instance $S_1,$ as shown in figure \ref{box},
the surface integrals 
\begin{equation}\label{ss1}
g_1=\iint_{S_1} \frac{1}{|{\bf x}'-{\bf x}_p|}\,\mathrm{d}S,
\end{equation}
\begin{equation}\label{ss2}
 {\bf g}_2=\iint_{S_1} \frac{{\bf x}'-{\bf x}_p}{|{\bf x}'-{\bf x}_p|^2}\,\mathrm{d}S,
\end{equation}
\begin{equation}\label{ss3}
 {\bf g}_3=\iint_{S_1} \frac{{\bf x}'-{\bf x}_p}{|{\bf x}'-{\bf x}_p|^3}\,\mathrm{d}S,
\end{equation}
are singular where
$${\bf x}'-{\bf x}_p=(0, y-y', z-z').$$ 
Defining 
$$ \bar {\bf r}=(y'-y_j, z'-z_k),$$
and  $$\bar r=|\bar {\bf r}|,$$
we apply the divergence theorem on (\ref{ss1}), thus we need to find a function $\varphi(\bar r)$ that satisfies 
$$\nabla \bar {\bf r}\varphi( \bar r))=\frac{1}{ \bar r},$$
or equivalently $$2\varphi( \bar r)+ \bar r\varphi'( \bar r)=\frac{1}{ \bar r}.$$ 
Solving the above equation, we get
 $$\varphi( \bar r)=\frac{1}{ \bar r}.$$
Thus $g_1$ is turned into 
\begin{equation*}
g_1=\sum_{n=1}^4\int_{L_n}\frac{{\bf x}'-{\bar {\bf x}}}{|{\bf x}'-{\bar {\bf x}}|} \cdot {\bf n}_n\,\mathrm{d} L,
\end{equation*}
where ${\bar {\bf x}}=(y_j, z_k),$ ${\bf n}_n$ is  the unit normal of $L_n.$  There is no singularity any more and $g_1$ can be calculated 
using numerical integration.

For (\ref{ss2}) and (\ref{ss3}),
due to the symmetry of vector ${\bf x}'-{\bf x}_p$ on $S_1,$ we have 
$$ {\bf g}_2=(0, 0, 0),$$ and 
$$ {\bf g}_3=(0, 0, 0).$$
\section{Parallelization}\label{Parallelization}
This paper closely follows \cite{Aihua2} and we therefore directly address the final numerical solving system of the EOS formulations of the 3D Maxwell's equations.
For the inside domain, the updating rule follows (\ref{3dmatrixequation}).
For the boundary part, the discretized boundary integral identities are represented by 
\begin{align}\label{boundarydiscretization3}
M_1  \left(
\begin{array}
[c]{c}%
{\bf E}_p^n \\
{\bf B}_p^n %
\end{array}
\right)&  =\left(
\begin{array}
[c]{c}%
{\bf E}_R \\
{\bf B}_R %
\end{array}
\right),
\end{align}
where 
$\left( \begin{array}
[c]{c}%
{\bf E}_p^n\\
{\bf B}_p^n %
\end{array} \right)$ are the solution at the surface point ${\bf x}_p$ at time $t^n$ and $\left( \begin{array}
[c]{c}%
{\bf E}_R \\
{\bf B}_R %
\end{array} \right)$ are the summations of the integrals in the boundary integral representations after moving the unknowns to the left of the equations.
From equation (\ref{3dmatrixequation}), it is easy to see that the updating for the inside domain at time now will only involve the values that are one time step before.
While the solutions on the surface point ${\bf x}_p$ in (\ref{boundarydiscretization3}) require both the historical values of the current density and the charge density and the historical field values of all the surface points due to the retarded integrals involved. Therefore the part of the code calculating the  surface solution dominates both the memory usage and the processor usage.
The calculations are therefore parallelized based on partitioning the surface into pieces and distributing each piece to separate processors, whereas the inside of the scattering object is residing on each processor.
The updating processes are illustrated by the following C code where
   \begin{equation*}
\begin{split}
&\text{p} : \quad \text{index of surface point}\\
&\text{n} : \quad \text{index of time level}\\
&\text{es, bs} : \quad \text{fields solutions on surface up to time}\, \,t^{n-1}\\
&\text{e, b} : \quad \text{fields solutions of inside domain at time}\, \,t^{n}\\
&\text{el, bl} : \quad \text{fields solutions of inside domain at time}\, \,t^{n-1}\\
&\text{J, P} : \quad \text{current density and electric density up to time}\, \,t^{n-1}\\
&\text{UpdateS(p,n,J,P,es,bs)}: \quad \text{update surface solutions at ${\bf x}_p$ at time} \, \,t^{n}\\
&\text{UpdateV(e,b,el,bl,es,bs,J,P,n)}: \quad \text{update inside solutions at time} \, \,t^{n}
\end{split}
  \end{equation*}
\begin{lstlisting}[style=CStyle]

int rank,size;// processor id and number of processors

MPI_Init(&argc,&argv);
MPI_Comm_size(MPI_COMM_WORLD,&size);
MPI_Comm_rank(MPI_COMM_WORLD,&rank);

int Nt,Ns,Nss,lp,lsize;
int p,index,indexeb;

lp=Ns/size;
lsize=3*lp;
double lse[lsize],lsb[lsize];	

double *esurface,*bsurface;
esurface = (double *)malloc(Nss*sizeof(double)); 
bsurface = (double *)malloc(Nss*sizeof(double)); 

for(n=0;n<Nt;n++){
  for(p=rank*lp;p<(rank+1)*lp;p++){
      gsl_vector *lresult=gsl_vector_alloc(6);
     // update the surface values at each grid point in parallel
      lresult=UpdateS(p,n,J,P,es,bs);
      index=(p%lp)*3;
      for(i=0;i<3;i++){
        indexeb=index+i;
        lse[indexeb]=gsl_vector_get(lresult,i);
        lsb[indexeb]=gsl_vector_get(lresult,i+3);
       
      }
     gsl_vector_free(lresult);
   }
   //collect data from all processes
    MPI_Allgather( lse, lsize, MPI_DOUBLE, esurface, lsize, MPI_DOUBLE, MPI_COMM_WORLD);
    MPI_Allgather( lsb, lsize, MPI_DOUBLE, bsurface, lsize, MPI_DOUBLE, MPI_COMM_WORLD);

   //updating the whole surface
    for(p=0;p<Nss;p++){
      gsl_matrix_set(es,p,n,*(esurface+p));
      gsl_matrix_set(bs,p,n,*(bsurface+p));
    }

  //update the inside domain by the domain based method supported by the surface values
  UpdateV(e,b,el,bl,es,bs,J,P,n); 

}

MPI_Finalize();

\end{lstlisting}

\end{appendices}

\bibliographystyle{unsrt}
\bibliography{paper3}

\end{document}